\documentclass[12pt]{article}

\usepackage{graphicx, amssymb, latexsym, amsfonts, amsmath, lscape, amscd,
amsthm, color, epsfig, mathrsfs, tikz, enumerate, multirow}
\usepackage{hyperref}
\usepackage{cleveref}
\usepackage{bm}
\usepackage{todonotes}
\usetikzlibrary{shapes.misc}
\tikzstyle{noeud}=[circle,inner sep=2, minimum size =3 pt, line width = 1pt, draw=black, fill=white]
\usepackage{subcaption}
\newtheorem{theorem}{Theorem}[section]

\newtheorem{corollary}[theorem]{Corollary}

\newtheorem{lemma}[theorem]{Lemma}
\newtheorem{proposition}[theorem]{Proposition}
\newtheorem{question}[theorem]{Question}

\newtheorem{definition}{Definition.}
\newtheorem{observation}[theorem]{Observation}

\newtheorem{remark}[theorem]{Remark}
\newtheorem{claim}{Claim}[theorem]
\newenvironment{answer}
  {\par\noindent\textbf{Answer.}\ }
  {\par}

\newcommand{\qedclaim}{\hfill $\diamond$ \medskip}

\newcommand\DELETE[1]{}

\title{On Monitoring Edge-Geodetic Sets of Dynamic Graph}
\author{
{\sc Zin Mar Myint}$\,^{(a,} \,\thanks{The first author was supported by the Doctoral Fellowship in India for ASEAN (DIA:2020-25).}\,\,^{)}$ , {\sc Ashish Saxena}$\,^{(b)}$\\
\mbox{}\\
{\small $(a)$ Indian Institute of Technology Dharwad, India.}\\
{\small $(b)$ Indian Institute of Technology Ropar, India.}\\}

\date{}

\begin{document}

\maketitle

\begin{abstract}
The concept of a monitoring edge-geodetic set (MEG-set) in a graph $G$, denoted $MEG(G)$, refers to a subset of vertices $MEG(G)\subseteq V(G)$ such that every edge $e$ in $G$ is monitored by some pair of vertices $ u, v \in MEG(G)$, where $e$ lies on all shortest paths between $u$ and $v$. The minimum number of vertices required to form such a set is called the monitoring edge-geodetic number, denoted $meg(G)$. 
The primary motivation for studying $MEG$-sets in previous works arises from scenarios in which certain edges are removed from $G$. In these cases, the vertices of the $MEG$-set are responsible for detecting these deletions. 
Such detection is crucial for identifying which edges have been removed from $G$ and need to be repaired. In real life, repairing these edges may be costly, or sometimes it is impossible to repair edges. In this case, the original $MEG$-set may no longer be effective in monitoring the modified graph. This highlights the importance of reassessing and adapting the $MEG$-set after edge deletions. 
This work investigates the monitoring edge-geodetic properties of graphs, focusing on how the removal of $k$ edges affects the structure of a graph and influences its monitoring capabilities. Specifically, we explore how the monitoring edge-geodetic number $meg(G)$ changes when $k$ edges are removed. The study aims to compare the monitoring properties of the original graph with those of the modified graph and to understand the impact of edge deletions.
\end{abstract}


\section{Introduction}
In a network, edges signify the communication links between vertices, and monitoring them is crucial for several reasons. First, edges serve as the pathways for data flow, and any failure or malfunction in these links can disrupt the entire network's performance. Monitoring edges help ensure reliable data transmission, detect bottlenecks, and prevent potential breakdowns in connectivity. Moreover, edge monitoring can identify security vulnerabilities such as unauthorized access or data breaches that may occur via compromised links.

In 2023, Foucaud et al. \cite{foucaud2023monitoring} introduced an innovative graph-theoretic concept known as ``monitoring edge-geodetic sets" (shortly, MEG-sets). This concept addresses the challenge of network monitoring, specifically focusing on the detection and repair of faults within a network when certain connections (edges) fail. In their paper \cite{foucaud2023monitoring}, the authors define $G$ as a finite, undirected, simple connected graph. Their objective is to identify the monitoring edge-geodetic set of graph $G$. To achieve this, we select a small subset of vertices (representing the probes) of the network such that all connections are covered by the shortest paths between pairs of vertices in the network.
Moreover, any two probes are able
to detect the current distance that separates them. The objective is for a pair of probes to detect a change in their distance value when an edge of the graph is deleted, thus allowing us to detect the failure. Specifically, for any two vertices $u$ and $v$ in $G$, an edge $e$ is considered monitored by $u$ and $v$ if it lies on all shortest paths between them. This is a well-studied problem \cite{foucaud2023monitoring, foucaud2023monitoringfull,foucaud2024bounds,Foucad_2024,haslegrave2023monitoring,xu2024monitoring,Ma_2024,Davide_2024, Li_2024}. The formal definition of MEG-sets of the graph $G$ is as follows.
\begin{definition}\cite{foucaud2023monitoring}
A subset of vertices of graph $G$, say $M$, is called a Monitoring Edge-Geodetic set (MEG-set) if, for every edge $e$ in $G$, there exists a pair of vertices in $M$ that monitors $e$. In other words, for any two vertices $u$ and $v$ in $M$, edge $e$ lies on all shortest paths between them. A MEG-set of $G$ is denoted by $MEG(G)$, and the minimum size of the $MEG(G)$ is denoted by $meg(G)$.
\end{definition}

Recently, Foucaud et al. \cite{Foucad_2024}, \cite{foucaud2024bounds} studied the subdivision operation, where each edge of a graph $G$ is subdivided $k$ times by adding the $k$ vertices to each edge of $G$ which is an example of a dynamic graph operation. Subdividing to the edges means adding the new vertices to the existing graph, which can have a direct effect on the monitoring edge-geodetic number.
As the graph evolves dynamically through these operations, the monitoring edge-geodetic sets must adapt to the new structure on $meg(G)$, providing insight into how the graph's monitoring properties change with each transformation. Specifically, as mentioned earlier, probes placed as part of a $MEG$ set can detect changes in shortest path distances, providing critical information on edge deletion. It might be possible the faulty edge can not be repaired. In this case, the initial placement of probes in $G$ may not be sufficient to monitor the new state of the network. This raises important questions: How many probes are required for effective monitoring after changes? Should the placement of these probes be adjusted in the updated network? Answering these questions leads us to investigate $MEG$-sets in dynamic graphs, where continuous adaptation and strategic probe placement are crucial. Before we study the problem of $MEG$-sets in dynamic networks, we need to understand the formal description of dynamic graphs.
\begin{definition}\cite{harary1997dynamic}
    A dynamic graph is a sequence of graphs $G_0, \, \ldots,\, G_m$ on $n$ vertices such that $G_{i+1}$ is obtained from $G_i$ by adding or removing edges or vertices.
\end{definition}
These graphs evolve as connections are established, broken, or rerouted due to factors like vertex mobility, edge failures, or the addition of new resources \cite{Jacob_1998, Camil_2004, Aaron_2011, Sankowski_2004, Sayan_2017}. Numerous dynamic graph problems, such as the problem of dynamic shortest paths, have been tackled in multiple studies: In 2004, Demetrescu et al. \cite{Camil_2004} proposed efficient algorithms for maintaining shortest path information in graphs that frequently change. Similarly, in 2011, Bernstein et al. \cite{Aaron_2011} focused on the development of methods for quickly updating shortest path calculations in response to dynamic changes in the graph. The vertex cover problem has also been a focal point in dynamic graph studies. In 2017, Bhattacharya et al. \cite{Sayan_2017} investigated algorithms for maintaining approximate vertex covers in dynamic graphs, providing solutions that balance efficiency with the need for accuracy in changing scenarios. Their work contributes to understanding how to cover vertices effectively while accommodating graph modifications.

In this paper, we investigate the Monitoring Edge-Geodetic problem on dynamic graphs. Specifically, we explore this problem on dynamic graphs where edges can be deleted. To the best of our knowledge, this problem has not been studied. 

\begin{remark}
     One important thing to observe if some edge(s) are removed from the graph, then the resultant graph may be disconnected. In this case, the $MEG$-set in the resultant graph is the union of the $MEG$-set of each connected component. If the component is an isolated vertex, then the $MEG$-set is considered to be an empty set because it does not contribute to monitoring the edges.
\end{remark}

\noindent \textbf{Related notions.} We denote the degree of a vertex $u$ by $deg(u)$. A MEG-set of $G$ is denoted by $MEG(G)$, and the minimum size of the $MEG(G)$ is denoted by $meg(G)$.

\subsection{Contributions and Structure of the Paper}
In this paper, we make several key contributions to the study of monitoring edge-geodetic sets and numbers of dynamic graphs, particularly focusing on the effects and properties of edge removal within graph structures. Our contributions and the structure of the paper are as follows:
\begin{itemize}
    \item In Section~\ref{preliminaries}, we recall some definitions and theorems from the previous research, which we reference while proving our results.
    \item In Section~\ref{dynamic}, for each of the graph classes mentioned (path, cycle, unicyclic graph, tree, grid graph, etc.), we explore the changes in the monitoring edge-geodetic number when an edge is removed. We provide specific results and insights into how edge removal alters the $meg$ in these graphs. Moreover, we extend the analysis to the scenario where multiple edges are removed from trees. Specifically, we investigate the boundedness of the $meg$ for the resultant graph $G$ after removing at most $k$ edges from a tree.
    \item In Section~\ref{general}, we explore the $meg$ of general graphs. Specifically, how the removal of cut edge edges and pendant edges (cut-edges with one endpoint having degree one) and an edge $e$ incident to vertex a simplicial vertex is removed from a connected graph $G$ impacts the monitoring edge-geodetic number. In this section, we provide the bounds on the $meg$ for graphs after the removal of these edges.
    \item In Section~\ref{sec:conclusion}, we conclude our work and discuss the future scope.

\end{itemize}

\section{Preliminaries results}\label{preliminaries}
Let $G$ be connected graph. A \emph{cut edge} in a graph $G$ is an edge whose removal increases the number of connected components of $G$, and a pendant edge $e$ in a graph $G$ is a cut edge with one of the end vertex of edge $e$ with degree one. A \emph{cut vertex} is a vertex in a graph that, when removed, increases the number of components in the graph. A \emph{simplicial vertex} in $G$ is a vertex such that its neighborhood induces a clique in $G$.

In 2023, Foucad et al. introduced the $MEG$-sets \cite{foucaud2023monitoring}; the authors determined the value of $meg(G)$ for several graph classes. In this work, we build upon several of their results, which are crucial for understanding the subsequent developments in our study. Some of these results are expressed as follows: 

\begin{lemma} \label{lm:existing-simplicial}
    \cite{foucaud2023monitoring}In a graph $G$ with at least one edge, any simplicial
vertex belongs to any edge-geodetic set and, thus, to any $MEG(G)$.
\end{lemma}
\begin{theorem}\label{thm:existing-pendant}
    \cite{foucaud2023monitoring} A pendant vertex is a part of every $MEG(G)$.
\end{theorem}

\begin{theorem}\label{thm:existing-Tree}
    \cite{foucaud2023monitoring} For any tree $T$ with at least one edge, the only optimal $MEG(T)$ consists of the set of leaves (or pendant vertices) of $T$.
\end{theorem}

\begin{corollary}\label{cor:existing-path}
    \cite{foucaud2023monitoring} For any path $P_n$ with at least one edge, $meg(P_n)=2$.
\end{corollary}

\begin{theorem}\label{thm:existing-Cycle}
    \cite{foucaud2023monitoring} Given an $n$-cycle graph $C_n$, for $n= 3$ and $n \geq 5$, $meg(C_n)=3$. Moreover, $meg(C_4)=4$.
\end{theorem}

\begin{theorem}\label{thm:unicycle}\cite{foucaud2023monitoring}
Let $G$ be a unicyclic graph where the only cycle $C^*$ has length $k$ and whose set of pendant
vertices is $L(G),\ |L(G)|=l$. Let $V^+_{c^*}$
be the set of vertices of $C^*$ with degree at least $3$. Let $p(G)=1$ if $G[V(C^*)\backslash V^+_{c^*}]$ contains a path whose length is at least $\lfloor \frac{k}{2}\rfloor$, and $p(G)=0$ otherwise. Then, if $k\in \{3, 4\},\
meg(G)=l+k-|V^+_{c^*}|$. Otherwise $(k\geq 5)$, then 
\begin{equation*}
meg(G) =
\begin{cases} 
 3, & \text{if $|V^+_{c^*}|=0$;},\\
 l+2, & \text{if $|V^+_{c^*}|=1$;},\\
  l+p(G)+1, & \text{if $|V^+_{c^*}|=2$, $k$ is even, and the vertices in $V^+_{c^*}$ are}\\
   & \text{adjacent or opposite on $C^*$;}\\
   l+p(G), & \text{in all other cases.}
   \end{cases}
\end{equation*}
\end{theorem}

\begin{theorem}\cite{foucaud2023monitoring}\label{rectangular}
Let $G$ be a rectangular grid of size $m\times n$. For any $m,n\geq 2$, then $meg (G)=2(m+n-2)$.
\end{theorem}

\begin{lemma}\label{lm:existing-cut_vertex}
    \cite{foucaud2023monitoring} Let $G$ be a graph, and $u$ be a cut-vertex of $G$. Then $u$
is never part of any minimal $MEG(G)$.
\end{lemma}

Those results have been examined in several studies \cite{foucaud2023monitoringfull,foucaud2024bounds,Foucad_2024,haslegrave2023monitoring,xu2024monitoring,Ma_2024,Davide_2024,Li_2024}. An important question raised in the conclusion of \cite{foucaud2023monitoring} is: which vertices are always included in every $MEG(G)$? In 2024, Foucad et al. \cite{Foucad_2024, foucaud2024bounds} addressed this question, and their result is as follows.

\begin{theorem}\label{thm meg full}
    \cite{Foucad_2024, foucaud2024bounds} Let $G$ be a graph. A vertex $v \in V (G)$ is in every $MEG(G)$ if
and only if there exists $u \in N(v)$ such that any induced $2$-path $uvx$ is part of a $4$-cycle.
\end{theorem}

 Additionally, the authors have explored various bounds on \(meg(G)\) and conducted their research across different graph classes. We will utilize some of their findings in our work, which are outlined as follows. A split graph G is a graph whose vertices can be partitioned into a clique and an independent set.

\begin{corollary}\label{cor:exiting-split}
    \cite{foucaud2024bounds} Let $G$ be a split graph with $k$ vertices having a pendent neighbour. If $G$ has $n$ vertices, then $meg(G) = n - k$.
\end{corollary}

\begin{proposition}\cite{foucaud2024bounds}\label{bridge}
 There exists the vertices $u, v\in V(G)$ such that $uv$ is a cut edge of the graph $G$, then any vertex $u,v$ is either a vertex having the degree
 $1$, or is never part of any minimum $MEG(G)$.
\end{proposition}


\section{Results on the various graphs classes}\label{dynamic}
In this section, we analyze the monitoring edge-geodetic number of a graph $G$, i.e., $meg(G)$, in dynamic graphs across different graph classes, including trees, paths, cycles, unicyclic graphs, grid graphs, etc. In a dynamic graph scenario, edges are added or removed over time, which can lead to significant changes in the graph structure. However, in this study, we focus specifically on the effect of edge removal within these particular graph classes. For instance, removing an edge from a cycle transforms it into a path, while removing an edge from a path results in a forest. Now, we will study the following $meg$ of the dynamic graph by removing some edges from the particular graphs.

 

\subsection{Trees}
Let $G$ be a tree. This section discusses how $meg$ changes if the edge(s) is removed from $G$. We start by considering the scenario when an edge is removed from $G$. If an edge $e=uv$ is removed from a tree $G$, we obtain the following for $meg$ of the graph $G \backslash \{e\}$, denoted $G'$.

\begin{theorem}\label{lm:tree-1-dyn}
    Let $G$ be the tree with order $n$. If we remove an edge $uv\in E(G)$ from the graph $G$, then $meg(G')$ is as follows. 
    \begin{itemize}
        \item[(i)] If the $deg(u)=deg(v)=1$, then $meg(G')=0$.  
        \item [(ii)] If either $deg(u)=1$ and $deg(v)=2$ or $deg(u)\geq 3$ and $deg(v)\geq 3$, then $meg(G')= meg(G)$.
        \item[(iii)] If $deg(u)=deg(v)=2$, then $meg(G')= meg(G)+2$.
        \item[(iv)] If we remove an edge which has $deg(u)=2$ and $deg(v)\geq 3$, then $meg(G')=meg(G)+1$.
        \item[(v)] If $deg(u)=1$ and $deg(v)\geq 3$, then $meg(G')=meg(G)-1$.
    \end{itemize}
\end{theorem}

\begin{proof}
  Due to Theorem \ref{thm:existing-Tree}, the $meg(G)$ corresponds to the number of pendant vertices in $G$. After removing an edge from the tree, the resulting graph $G'$
  becomes a forest. It’s important to note that after removing edge $e$, the graph $G'$ may consist of at most two trees. Therefore, $meg(G')$ will be determined by the number of pendant vertices in $G'$.  

\noindent \textit{Case (1):} If $deg(u)=deg(v)=1$, then $G$ is a path of length 1. If the edge $uv$ is deleted, then there is no edge to monitor. Therefore, $meg(G')=0$.

\noindent \textit{Case (2):} If $deg(u)=deg
        (v)\geq 3$, the number of pendant vertices remains unchanged after removing $e$. Thus  $meg(G')=meg(G)$.

\noindent \textit{Case (3):} If $deg(u)=2$ and $deg(v) \geq 3$, then after removing edge $e$, the number of pendant vertices increases by 1. Therefore, $meg(G')=meg(G)+1$. Similarly, in the case, when $deg(u)\geq 3$ and $deg(v) =2$, one can prove $meg(G')=meg(G)+1$.

\noindent \textit{Case (4):} If $deg(u)=deg(v)=2$, then after removing edge $e$ from $G$, the number of pendant vertices increases by 2. Therefore, $meg(G')=meg(G)+2$.

\noindent \textit{Case (5):} If $deg(u)=2$ and $deg(v)=1$, then after removing edge $e$ from $G$, vertex $v$ becomes an isolated vertex and vertex $u$ becomes an pendant vertex. Therefore, the number of pendant vertices remains the same, and $meg(G')=meg(G)$. Similarly, in the case, when $deg(u)=1,$ and $deg(v)=2$, one can prove $meg(G')=meg(G)$. 

\noindent \textit{Case (6):} If $deg(u)\geq 3$ and $deg(v)=1$, then after removing edge $e$ from $G$, vertex $v$ becomes an isolated vertex and vertex $u$ does not become an pendant vertex. Therefore, the number of pendant vertices decreases by 1, and $meg(G')=meg(G)-1$. Similarly, in the case, when $deg(u)=1 $ and $deg(v)\geq 3$, one can prove $meg(G')=meg(G)-1$. 
    
    This completes the proof. 
\end{proof}

Based on the result on trees, we have the following result on a path of length $n$, say $P_n$, and a cycle of length $n$, say $C_n$. Suppose an edge from $P_n$ and $C_n$ is removed, and the resultant graph is $P'_n$ and $C'_n$, respectively. It is important to note that $P'_n$ is a forest (nothing but a collection of trees), and $C_n'$ is a tree. Due to Corollary \ref{cor:existing-path}, we know that $meg(P_n)$ is $2$. And, due to Theorem \ref{thm:existing-Cycle}, we know that $meg(C_n)$ is $3$ or $4$ based on the value of $n$.

\begin{corollary}\label{cor:path}
    After removing an edge $uv$ from $P_n$, we have the following results such that 
    \begin{itemize}
        \item[(a)] If $deg(u)=deg(v)=1$, then $meg(P'_n)=0$.
        \item[(b)] If $deg(u)=1$ and $deg(v)=2$, then $meg(P'_n)=meg(P_n)=2$.
        \item[(c)] If $deg(u)=deg(v)=2$, then $meg(P'_n)=meg(P_n)+2=4$.
    \end{itemize}
\end{corollary}
\begin{proof}

Since $P_n$ is a tree, therefore it is an implication of Theorem \ref{lm:tree-1-dyn}.
\end{proof}

\begin{corollary}\label{cor:cycle}
    After removing an edge $uv$ from $C_n$, $meg(C'_n)=2$.
\end{corollary}
\begin{proof}
    After removing an edge from $C_n$, $C'_n$ is nothing but a path of length at least 2. And, we know $meg$ of a path is $2$ (using Corollary \ref{cor:existing-path}). Therefore, $meg(C'_n)=2$.  
\end{proof}

Now we consider the based on the case when at most $k$ edges are removed from $G$. Let $G'$ be the resultant graph after removing $k$ edges from $G$. We have the following result. 

\begin{theorem}\label{full_edges_removal}
Let $G$ be a tree, and let $G'$ be the resultant graph after removing at most $k$ edges from $G$. Then, the following inequality holds:
$$0 \leq meg(G') \leq meg(G) + 2k.$$
\end{theorem}

\begin{proof}
     Due to Theorem \ref{thm:existing-Tree}, we know that the $meg(G)$ is the number of pendant vertices in $G$. We know that after removing edge(s) from $G$, it becomes a forest (nothing but a collection of trees). Therefore, $meg(G')$ is nothing but the number of pendant vertices. The number of edges in $G$ is $|V|-1$. If $k \geq |V|-1$, then $G'$ contains $|V|$ isolated vertices. Therefore, $meg(G')\geq 0$. In Theorem~\ref{lm:tree-1-dyn}, 
      we established that the removal of each edge can create at most $2$ additional pendant vertices. Therefore, after removing $k$ edges, the maximum increase in pendant vertices is $2k$. Thus, we have: $meg(G') \leq meg(G)+2k$. Combining these results completes the proof.
\end{proof}
Building on this result, Corollary~\ref{half_edge_removal} provides a more refined bound result as below.

\begin{corollary} \label{half_edge_removal}
    If $k > \lceil \frac{|E|}{2}\rceil-1$, then $meg(G')<meg(G)+2k$.  
\end{corollary}
\begin{proof}
     Suppose the given statement is false. Therefore, $meg(G')=meg(G)+2k$ for $k > \lceil \frac{|E|}{2}\rceil-1$. Let edges $e_1$, $e_2$, \ldots, $e_k$ be removed from $G$, and $e_i=(u_i, \, v_i)$. As per Theorem~\ref{lm:tree-1-dyn}, removing an edge $(u, \,v)$ increases the number of pendant vertices by $2$, if $deg(u)=deg(v)=2$. Hence, for edge $e_i$, $deg(u_i)=2$ and $deg(v_i)=2$. If for $i \neq j$, $e_i$ and $e_j$ are incident edges, then the number of pendant vertices does not increase by $4$. It is because removing edge $e_i$ and $e_j$ creates an isolated vertex. Therefore, $deg(u_i)=deg(v_i)=2$, and edges $e_i$ and $e_j$ are not incident for any $i \neq j$. In this case, the number of vertices in $G'$ is at least $|V|+1$, which is a contradiction due to the following reasons. In $G$, there are at least $2k+2$ vertices present due to $meg(G')=meg(G)+2k$, and $meg(G)\geq2$. Without loss of generality, let $k=\lceil \frac{|E|}{2}\rceil$. In this case, the number of vertices are $2\times\lceil \frac{|E|}{2}\rceil+2=2 \times \lceil \frac{|V|-1}{2}\rceil+2\geq |V|+1$. Therefore, our assumption is wrong. Hence, $meg(G')<meg(G)+2k.$ 
\end{proof}


The next question we would like to pose is as follows: for every integer \( k \) where \( 0 \leq k \leq \lceil \frac{|E|}{2} \rceil - 1 \), can we identify a class of trees \( G \) such that \( meg(G') = meg(G) + 2k \) after removing the edges \( e_1, e_2, \ldots, e_k \)? In other words, removing \( k \) edges should increase the \( meg \) value by \( 2k \). The construction for such a tree \( G \) is described as follows.

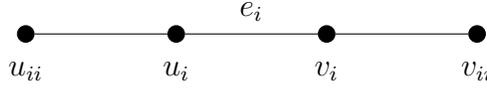
\begin{figure}[h!]
			\centering
			\begin{subfigure}[b]{\textwidth}
				\centering
					\begin{tikzpicture}
        \node[noeud,fill=black,circle] (a) at (0,0){};
        \node[circle] (u_ii) at (0,-0.5) {$u_{ii}$};
        
        \node[noeud,fill=black,circle] (b) at (2,0){};
        \node[circle] (u_i) at (2,-0.5) {$u_i$};
        
        \node[noeud,fill=black,circle] (c) at (4,0){};
        \node[circle] (v_i) at (4,-0.5) {$v_i$};

        \node[noeud,fill=black,circle] (d) at (6,0){};
        \node[circle] (v_ii) at (6,-0.5) {$v_{ii}$};
    
        \draw (a) -- (b) -- (c) -- (d);
    
        \node at (3,0.3) {$e_i$};
    \end{tikzpicture}
\caption{Block $B_i$.}
\label{fig:block}
\end{subfigure}	\hfill	
\begin{subfigure}[b]{\textwidth}
\centering
	\begin{tikzpicture}[scale=0.5, transform shape]
\node at (-1.7,1.5){\huge$T_1$};
\node at (3.2,1.5){\huge$T_2$};
\node at (12.2,1.5){\huge$T_k$};
\node at (17.2,1.5){\huge$T_{k+1}$};

\node at (2.5,-0.7){\huge$e_1$};
\node at (16.75,-0.7){\huge$e_k$};
\draw[fill=none, dotted, very thick] (0.5,0) circle (0.3) node {};

\node[noeud,fill=black,circle] (a1) at (0.5,0){};
\node[noeud,fill=black,circle] (a11) at (2,0){};
\node[noeud,fill=black,circle] (a12) at (3,0){};
\node[noeud,fill=black,circle] (a13) at (4.5,0){};

\draw[fill=none, dotted, very thick] (4.5,0) circle (0.3) node {};

\draw[fill=none, dotted, very thick] (5.5,0) circle (0.3) node {};

\node[noeud,fill=black,circle] (b1) at (5.5,0){};
\node[noeud,fill=black,circle] (b2) at (4.5,0){};
\node[noeud,fill=black,circle] (b11) at (7,0){};
\node[noeud,fill=black,circle] (b12) at (8,0){};



\node[noeud,fill=black,circle] (b13) at (11,0){};
\node[noeud,fill=black,circle] (b14) at (12,0){};
\draw[fill=none, dotted, very thick] (13.5,0) circle (0.3) node {};
\draw[fill=none, dotted, very thick] (14.7,0) circle (0.3) node {};

\node[noeud,fill=black,circle] (c1) at (13.5,0){};
\node[noeud,fill=black,circle] (c11) at (16,0){};
\node[noeud,fill=black,circle] (c12) at (17,0){};

\node[noeud,fill=black,circle] (c2) at (14.7,0){};
\node[noeud,fill=black,circle] (d1) at (18.5,0){};

\draw[fill=none, dotted, very thick] (18.5,0) circle (0.3) node {};
\node (a) at (0,0) {\begin{tikzpicture}[scale=0.6, transform shape]
\draw plot [smooth cycle, tension=1] coordinates { (-4.75-.5,.5+.05) (-4.7-.1,-1-.5) (-2+.2,-1.5-.5) (-1.5+.5,.5) (-2.8+.5,1.5+.5) (-4.4,1.5+.5)};
\node at (-2,-0.9){\huge $u_{11}$};

\end{tikzpicture}};
\node (b) at (5,0) {\begin{tikzpicture}[scale=0.6, transform shape]
\draw plot [smooth cycle, tension=1] coordinates { (-4.75-.5,.5+.05) (-4.7-.1,-1-.5) (-2+.2,-1.5-.5) (-1.5+.5,.5) (-2.8+.5,1.5+.5) (-4.4,1.5+.5)};
\node at (-4,-0.9){\huge $v_{11}$};
\node at (-2,-0.9){\huge $u_{22}$};
\end{tikzpicture}};
\node (c) at (14,0) {\begin{tikzpicture}[scale=0.6, transform shape]
\draw plot [smooth cycle, tension=1] coordinates { (-4.75-.5,.5+.05) (-4.7-.1,-1-.5) (-2+.2,-1.5-.5) (-1.5+.5,.5) (-2.8+.5,1.5+.5) (-4.4,1.5+.5)};
\node at (-3.9,-0.9){\huge$ v_{(k-1)(k-1)}$};
\node at (-1.5,-0.9){\huge $u_{kk}$};
\end{tikzpicture}};
\node (d) at (19,0) {\begin{tikzpicture}[scale=0.6, transform shape]
\draw plot [smooth cycle, tension=1] coordinates { (-4.75-.5,.5+.05) (-4.7-.1,-1-.5) (-2+.2,-1.5-.5) (-1.5+.5,.5) (-2.8+.5,1.5+.5) (-4.4,1.5+.5)};
\node at (-4,-0.9){\huge $v_{kk}$};

\end{tikzpicture}};

\draw (a1) -- (b2);
\draw (b1) -- (b12);
\draw[dashed] (b12) -- (b13);
\draw (b13) -- (c1);
\draw (c2) -- (d1);
\end{tikzpicture}
\caption{A vertex with dashed circle denotes the vertice $u_{ii}$ and $v_{ii}$, $i \in [1, \,k]$, $T_1$ , $T_{k+1}$ is a sub-tree of $T$ which contains vertex $u_{11}$, $v_{kk}$, respectively. A sub-tree $T_i$ contains vertex $u_{i-1 i-1}$ and $v_{i i}$, $i\in [2, k]$. For some $T_i$, $u_{i-1 i-1}$ and $v_{i i}$ can be same vertex, i.e., $T_i$ contains one vertex.}
\label{fig:tree-k-dyn}
\end{subfigure}
			\caption{(a) It is a block of graph $G$, (b) This graph denotes a tree where after removal edges $e_i$s, $i\in[1,\,k]$, the $meg$ increases by $2k$.}
\end{figure}

Suppose the edges $e_1$, $e_2$, \ldots, $e_k$ are removed from $G$, and $e_i=(u_i, \, v_i)$. Using Corollary~\ref{half_edge_removal}, edge $e_i$ should satisfy the following necessary conditions, which are as follows. For any $i \neq j$, $deg(u_i)=deg(v_i)=2$ and $e_i$ and $e_j$ are not incident edges. For edge $e_i$, we denote a block by $B_i$ (refer Fig. \ref{fig:block}). Therefore, $B_1$, $B_2$, \ldots, $B_k$ are blocks corresponding edges $e_1$, $e_2$, \ldots, $e_k$, respectively. For $i$, there is a sub-tree $T_{i+1}$ between the vertices $v_{ii}$ and $u_{i+1 \,i+1}$. Note that it is possible that vertices $v_{ii}$ and $u_{i+1 i+1}$ are the same vertices. We can see the tree $G$ in Fig. \ref{fig:tree-k-dyn}. In Fig. \ref{fig:tree-k-dyn}, if edges $e_1$, $e_2$, \ldots, $e_k$ are removed, then $meg(G')=meg(G)+2k$.

\subsection{Unicyclic graphs}

In this section, we provide the $meg$ of the resultant graph $G'$ after removing an edge from the unicyclic graph $G$.

\begin{corollary}
Let $G$ be a unicyclic graph where the only cycle $C^*$ has length $k$ and whose set of pendant
vertices is $L(G),\ |L(G)|=l$. Let $V^+_{c^*}$
be the set of vertices of $C^*$ with degree at least $3$. If $G'$ is obtained by removing an edge $e=uv$ from $C^*$, then  
\begin{equation*}
meg(G') =
\begin{cases} 
 l, & \text{if $u,v \in V^+_{c^*}$},\\
 l+1, & \text{if $u \notin V^+_{c^*}$ and $v\in V^+_{c^*}$ or (if $v \notin V^+_{c^*}$ and $u\in V^+_{c^*}$)},\\
 l+2, & \text{if $u,\,v \notin V^+_{c^*}$},\\
   \end{cases}
\end{equation*}      
\end{corollary}
\begin{proof}
    After the removal of the edge $e$, $G'$ be a tree. Using Theorem \ref{lm:tree-1-dyn}, $meg(G')$ is the number of pendant vertices. If $u, v \in V^+_{c^*}$, then the removal of edge $e$ does not increase the pendant vertices. Therefore, if $u, v \in V^+_{c^*}$, then $meg(G')=l$. If $u \notin V^+_{c^*}$ and $v\in V^+_{c^*}$ (or if $v \notin V^+_{c^*}$ and $u\in V^+_{c^*}$), then removal of edge $e$ increases one pendant vertex. Therefore, in this case, $meg(G')=l+1$. If $u,\,v \notin V^+_{c^*}$, then the removal of edge $e$ increases two pendant vertices. Therefore, $meg(G')=l+2$. This completes the proof. 
\end{proof}

\begin{observation}
    Let $G$ be a unicyclic graph where the only cycle $C^*$ has length $k$. If $G'$ is obtained by removing an edge $e=uv$, which is not part of $C^*$, then two components will exist. Let $G_1$ and $G_2$ be two components of $G$ after removal of edge $e$. It is important to note that either $G_1$ or $G_2$ is a unicyclic graph, and the other is a tree. Without loss of generality, let $G_1$ be an unicyclic graph and $G_2$ be a tree. Using Theorem \ref{thm:unicycle}, we can compute $meg(G_1)$, and using Theorem \ref{thm:existing-Tree}, we can compute $meg(G_2)$. Therefore, $meg(G')=meg(G_1)+meg(G_2)$.
\end{observation}



\subsection{Grids}
A \textit{rectangular grid} \( G \) is a finite, undirected graph consisting of vertices and edges arranged in an \( m \times n \) grid, where \( m \geq 2\) and \( n \geq 2\) are positive integers representing the number of rows and columns, respectively. Each vertex \( v \in G \) is denoted by its position \( (i, j) =(\) row index, column index$)$. In grid \( G \), the vertices are categorized with \textit{corner vertices} (\( C \)) are those located at the four corners of the grid, specifically \( (1, 1) \), \( ( 1, \,n) \), \( (m, 1) \), and \( ( m, \,n) \), each having degree $2$, \textit{edge vertices} (\( \mathcal{E} \)) are those located on the boundary of the grid but not at the corners. These vertices include those at positions \( ( 1, \,j) \) for \( 2 \leq j \leq n-1 \), \( ( m,\,j) \) for \( 2 \leq j \leq n-1 \), \( ( i, 1) \) for \( 2 \leq i \leq m-1 \), and \( ( i, n) \) for \( 2 \leq i \leq m-1 \), each having degree $3$, and \textit{internal vertices} (\( \mathcal{I} \)) are not on the boundary of the grid, located at positions \( (i, j) \) where \( 2 \leq i \leq m-1 \) and \( 2 \leq j \leq n-1 \). Each internal vertex has degree $4$.

   In Theorem \ref{rectangular}, it is stated that \( meg(G) = 2(m + n - 2) \) for \( m, n \geq 2 \). A graph is a partial grid if it is an arbitrary subgraph of a grid, not necessarily induced. In \cite{foucaud2023monitoring}, the authors studied $meg(G)$ of $G$. This work focuses on $meg(G')$, where \( G' \) is the partial grid obtained by removing an edge from \( G \).\\

 \noindent \textbf{Notation for grid:} Let $S=\{(i,j)\in V(G): \,i\in \{1,m\}\ and\ 1\leq j\leq n\ or\ j\in \{1,n\}\ and\ 1\leq j\leq m\}$, i.e., $S$ contains $\mathcal{C}$ and $\mathcal{E}$. For any edge $uv \in E(G)$, let $G'$ a partial grid obtained by removing an edge $uv$ from $G$. Let $N(u)$ represent the neighbourhood of vertex $u$. We define $c$ as the number of corner vertices in $N(u)$, i.e., $c = | N(u) \cap C |$. We denote $deg(u)$ and $deg(v)$ by $d_1$ and $d_2$, respectively. The results of $G'$ are as follows.

\begin{lemma}
    Let $G'$ be a partial grid obtained by removing an edge $uv$ from the rectangular $m \times n$ grid $G$, $m, \,n\geq2$. The following holds if $u$, $v\in \mathcal{I}$.

      \begin{equation*}
        meg(G') =
    \begin{cases} 
    k+|S_1|-2, & \text{if $S_1 \neq \emptyset$}\\
    k+2, & \text{if $S_1 =\emptyset$}\\
   \end{cases}
\end{equation*}
 Here, the set $S_1$ is defined as follows. Let $S_e=\{xy\in E(G): x \in N(u), \, y\in N(v)\; with\; deg(x)=deg(y)=3\}$. If $e=xy\in S_e$, then $x, \,y \in S_1$.
\end{lemma}

\begin{proof}
Let $meg(G)=k$. In \cite{foucaud2023monitoring}, authors have shown that $S$ of $2(m + n - 2)$ vertices of $G$ that form the boundary vertices of the grid form the only optimal $MEG(G)$. Hence, $k=2(m+n-2)$. Let $S'$=$S-S_1$. Before proving the lemma, we prove two claims which help to prove the lemma. 

    \begin{claim}\label{boundary_vertices}
        All vertices in the set $S'$ are still required in every $MEG(G')$.  
    \end{claim}
    
\begin{proof}[Proof of Claim]
Let $(i,j)\in S'$. We use Theorem \ref{thm meg full} to show that $(i,\,j)$ is in every $MEG(G')$. If $(i,\,j)\in \mathcal{C}$, without loss of generality, let $u_1=(1,1)$ be a corner vertex. The vertex $u_1$ has two neighbours $(2,1)$ and $(1,\, 2)$. Note that vertices $(1, 1)$, $(1,2)$, $(2,2)$, $(2,1)$ forms a 4 cycle in $G'$. Let $v_1=(1,2)$. In this case, $x=(2,1)$. Since induced $2$-path $v_1u_1x$ is part of a $4$-cycle, therefore, $(1,\,1)$ is part of every $MEG(G')$.

If $(i,j)\in \mathcal{E}$, without loss of generality, let $u_2=(i, 1)\in S'$ with $2\leq i\leq m-1$. Vertices $(i-1,1)$, $(i+1, 1)$, $(i, 1)$ are neighbours of $(i,1)$ in $G'$. Since $u_1\in S'$, $(i,1)$, $(i,2)$, $(i-1,2)$, $(i-1, 1)$ and $( i,1)$, $(i,2)$, $(i+1,2)$, $( i+1,1)$ form a 4 cycle. Let $v_2=(i,2)$. In this case, $x \in \{(i-1, 1), \,(i+1,1)\}$ such that if $x=(i-1, 1)$, then induced 2-path $v_2u_2x$ is part of 4-cycle due to the fact that $(i,1)$, $(i,2)$, $(i-1,2)$, $(i-1, 1)$ forms a 4-cycle and similarly, if $x=(i+1, 1)$, then induced 2-path $v_2u_2x$ is part of 4-cycle due to the fact that $( i,1)$, $(i,2)$, $(i+1,2)$, $( i+1,1)$ forms a 4-cycle. Therefore, $(i, 1)$ is in every $MEG(G')$. This completes the proof of our claim.
\end{proof}

\begin{claim}
        Vertices $u$ and $v$ are in every $MEG(G')$.  
\end{claim}
\begin{proof}[Proof of Claim]
    Without loss of generality, let $u=(i,j)$ and $v=(i+1, j)$. In $G'$, vertices of each set $M_1$, $M_2$, $M_3$ and $M_4$ forms a 4-cycle, where $M_1=\{(i-1,j-1),\,(i, j-1), (i,j), (i-1,j)\}$, $M_2=\{(i-1,j), \, (i,j),\, (i,j+1), (i-1, j+1)\}$, $M_3=\{(i+1,j-1), (i+2,j-1), (i+2,j), (i+1,j)\}$, $M_4=\{(i+1,j), (i+2,j), (i+2,j+1), (i+1,j+1)\}$. In graph $G'$, $N(u)=\{(i-1,\,j),  \,(i,\,j-1),\, (i,\,j+1)\}$, and $N(v)=\{ (i+1,j-1), \,(i+2,j)\, (i+1,j+1)\}$.
    
    \noindent$\rightarrow$ To apply Theorem \ref{thm meg full}, for vertex $u$, let $v_1=(i-1,j)$. In this case, $x=(i,\,j-1)$ and $(i,\,j+1)$ such that if $x=(i,j-1)$, then $2$--induced path $v_1ux$ is part of $4$-cycle due to $u,v_1,x \in M_1$ and similarly, if $x=(i,\,j+1)$, then $2$-induced path $v_1ux$ is part of $4$-cycle due to $u,v_1,x \in M_2$. 
    
    \noindent$\rightarrow$ To apply Theorem \ref{thm meg full}, for vertex $v$, let $v_2=(i+2,j)$. In this case, $x=(i+1,\,j-1)$ and $(i+1, \,j+1)$ sucht that if $x=(i+1,j-1)$, then $2$-induced path $v_2vx$ is part of $4$-cycle due to $v,v_2,x \in M_3$ and similarly, if $x=(i+1,\,j+1)$, then $2$-induced path $v_2vx$ is part of $4$-cycle due to $v,v_2,x \in M_4$. Hence, $u$ and $v$ are in every $MEG(G')$. 
\end{proof}

We divide the proof of lemma into two cases as follows.

\noindent \textit{Case (1)} ${S_1 =\emptyset}$: In this case, $S'=S$, therefore $S\cup \{u, \,v\}$ is in every $MEG(G')$. Without loss of generality, let $u=(i,j)$, $v=(i,j+1)$. Note that since $S_1 =\emptyset$, $3\leq i\leq m-2$. In $G'$, let us define the shortest path $P_1$, $P_2$, $P_3$, $P_4$, $P_5$, $P_6$, $P_7^l$, $P_8^{l_1}$ such that
    $P_1$ is shortest path between $(1,1)$ and $ (1, n)$, $P_2$ is shortest path between $(1,n)$ and $(m,n)$, $P_3$ is shortest path between $(m,n)$ and $(m, 1)$, $P_4$ is shortest path between $(m,1)$ to $(1,1)$, $P_5$ is shortest path between $(i,1)$ and $(i,j)$, $P_6$ is shortest path between $(i,j+1)$ and $(i, n)$, $P_7^l$ is shortest path between $(1,l)$ and $(m,l)$, $1\leq l \leq n$, and $P_8^{l_1}$ is shortest path between $(l_1, 1)$ and $(l_1, n)$, $1< l_1< m, \, l\neq i$. It is important to note that these shortest paths are unique in $G'$, and the end vertices of each path are in $S\cup\{u,\,v\}$. Each edge $e_1\in G'$ lies on at least one of these unique shortest paths. Since $S\cup \{u, \,v\}$ is in every $MEG(G')$, and each edge in $G'$ is monitored by $S\cup \{u,\,v\}$, therefore, $S\cup \{u, \,v\}$ is the smallest $MEG(G')$.

\noindent \textit{Case (2)} ${S_1\neq \emptyset}$: In this case, there is an edge $u_1v_1\in S_e$. Note that $deg(u_1)$ and $deg(v_1)$ is 3, and $u_1,v_1\notin S'$. We show that $S' \cup \{u, \,v\}$ is the smallest $MEG(G')$. Let $u_{11}$ be the neighbour of $u_1$, and $v_{11}$ be the neighbour of $v_1$. It is important to note that $u_{11}$,  $v_{11} \in S'$. All the incident edges of $u_1$ and $v_1$ are monitored due to the following fact. The incident edges of $u_1$ is $u_1u$, $u_1v_1$ and $u_{11} u_1$, and the incident edges of $v_1$ is $v_1v$, $u_1v_1$ and $v_{11} v_1$. There exists only one shortest path $P'_1=u_{11} u_1  v_1  v$ between the vertices $u_{11}$ and $v$, and the shortest path $P'_2=v_{11} v_1  u_1  u$  between the vertices $v_{11}$ and $u$. Therefor, path $P'_1$, and $P'_2$ monitors all the incident edges of $u_1$, $v_1$ due to $v_{11}, u_{11}$, $u$, $v$ is a part of every $MEG(G')$. All the remaining edges are monitored by vertices $S' \cup \{u, \,v\}$ due to the following fact. Like the earlier case, we can find the unique shortest paths between two vertices from set $S'\cup \{u, \,v\}$ such that all edges from $G'$ except neighbours of $u_1$ and $v_1$ can be monitored by these unique shortest paths. Therefore, $S' \cup \{u, \,v\}$ is the smallest $MEG(G')$, and $meg(G')=meg(G)+|S_1|-2$. Note that since there can be at most two edges in the set $S_e$ in $G$, $|S_1|=2$ or $4$. 

This completes the proof. 
\end{proof}



\begin{lemma}\label{grid_1}
   Let $G$ be a rectangular $m \times n$ grid graph with $meg(G) = k$, $m,n \geq 2$. Let $G'$ be the graph obtained by removing an edge $uv$ from $G$. The following is true.
     \begin{equation*}
meg(G') =
\begin{cases} 
 k-c, & \text{if $d_1 = 2$, $d_2 = 2$ (or if $d_3=3$ $d_2=4$ and $c\neq 0$}),\\
k-1, & \text{if $d_3=3,\, d_2=4$ and $c=0$}\\
   \end{cases}
\end{equation*}
where $c$ is $|N(u)\cap C|$ or $|N(v)\cap C|$.
\end{lemma}

\begin{proof}
We prove the above statement in two parts as follows.

\noindent \textit{Case (1) $ d_1=d_2=2:$} The rectangular grid $G$ must be either $2\times n$ or $m\times 2$. There are two sub-cases as follows.

    \begin{itemize}
        \item \textit{Sub-Case (1) ${G}$ is ${2\times 2}$ grid:} In this case, $|N(u)\cap C|=1$ and $|N(v)\cap C|=1$, i.e., $|N(u)\cap C|+|N(v)\cap C|=c=2$. Now we can observe that the vertices belong to $N(u)\cap C$ and the vertices belong to $N(v)\cap C$ lie on the shortest path between $u$ and $v$ in the resultant graph $G'$. Therefore, $meg(G')=meg(G)-|N(u)\cap C|-|N(v)\cap C|=k-c$.
        
        \item  \textit{Sub-Case (2) ${G}$ is ${2\times n}$ or ${m\times 2}$ grid, ${m,\,n\geq 3:}$} Let $u_1v_1$ be the parallel edge to $uv$. After removing edge $uv$, vertices $u$ and $v$ become pendant vertices. Using Theorem \ref{thm:pendant}, $u$ and $v$ are part of every $MEG(G')$. Using Claim \ref{boundary_vertices}, we can prove that set $S-\{u, \, v,\,u_1,\,v_1\}$ is part of every $MEG(G')$. Therefore, set $S'=S-\{u_1, \,v_1\}$ is part of every $MEG(G')$. We show $S'$ is the smallest $MEG(G')$ as follows. Without loss of generality, let $G$ be $2 \times n$ grid. In this case, $u=(1,1)$, $v=(2,1)$, $u_1=(1, 2)$ and $v_1=(2,2)$. Let us define the shortest path $P_1, P_2, P_3 P_4^l$ such that $P_1$ is the shortest path between $(1, 1)$ and $(1, n)$, $P_2$ is the shortest path between $(2, 1)$ and $(2, n)$, $P_3=(1,1) (1,2)(2,2)(2,1)$ is the shortest path between $(1,1)$ and $(2,1)$, and $P_4^l$ is the shortest path between $(1, l)$ and $(2, l)$, $3 \leq l \leq n$. It is important to note that these shortest paths are unique in $G'$, and the end vertices of each path are in $S'$. Each edge $e_1\in G'$ lies on at least one of these unique shortest paths. Since $S'$ is in every $MEG(G')$, and each edge in $G'$ is monitored by $S'$, therefore, the set $S'$ is the smallest $MEG(G')$, and $meg(G')=meg(G)-|N(u)\cap C|-|N(v)\cap C|=k-c$.
    \end{itemize}

\noindent \textit{Case (2) $d_1=3, d_2=4$:} Without loss of generality, let $d_1=3, d_2=4$. Let $N(u)=\{u_1,u_2,v\}$ and $N(v)=\{u,v_1,v_2,v_3\}$. We divide this case into three sub-cases.
    \begin{itemize}
         \item \textit{Sub-Case (1)} ${|N(u)\cap C|=2:}$ In this case, $G$ is either $3 \times n$ or $m \times 3$. Without loss of generality, let $G$ be $3 \times n$ grid. In this case, $u=(2,1)$ and $v=(2,2)$. Using Claim \ref{boundary_vertices}, we can show $S-\{(1,1),\,(2,1),\,(3,1),\,(1,2),\,(3,2) \}$ is part of every $MEG(G')$. Similarly, one can show that $v$ is in every $MEG(G')$. We show that $S_1=(S\cup \{v\})-\{(2,1),\,(1,2),\,(3,2)\}$ is the smallest $MEG(G')$. Set $S_1$ is forming the $MEG(G')$ due to the following reason. Let us define the shortest path $P_1$, $P_2$, $P_3$, $P_4$, $P_5$, $P_6^l$ such that $P_1$ is shortest path between $(1,1)$ and $ (1, n)$, $P_2$ is shortest path between $(3,1)$ and $(3,n)$, $P_3$ is shortest path between $(2,2)$ and $(2, n)$, $P_4$ is shortest path between $(3,1)$ to $(2,2)$, $P_5$ is shortest path between $(1,1)$ to $(2,2)$, and $P_6^l$ is shortest path between $(1,l)$ and $(3, l)$ for $1\leq l \leq n$, and $l\neq 2$. It is important to note that these shortest paths are unique in $G'$, and the end vertices of each path are in $S_1$. Each edge $e_1\in G'$ lies on at least one of these unique shortest paths. Therefore, set $S_1$ is a $MEG(G')$. Now, we will show $S_1$ is the smallest by showing that at least two vertices are needed in set $(S\cup \{v\}\})-\{(1,1),\,(2,1),\,(3,1),\,(1,2),\,(3,2) \}$ to form an $MEG(G')$.
        \vspace{0.2cm}

        Note that set $S_2=(S\cup \{v\})-\{(1,1),\,(2,1),\,(3,1),\,(1,2),\,(3,2) \}$ is not a $MEG(G')$ because edge $e$ between $(1,\,1)$ and $(2,\,1)$ is not monitored by set $S_2$. Therefore, we need at least one vertex in set $S_2$. Therefore, if $(S_2\cup \{ x\})$ forms an $MEG(G')$, then $x \in V(G')-S_2$ (all vertices from $\mathcal{I}$ with $\{(1,1),\,(2,1),\,(3,1),\,(1,2),\,(3,2) \}$). 

\vspace{0.2cm}
        $\rightarrow$ If $x \in \mathcal{I}$, then the edge $e_0$ between $(1,1)$ and $(2, 1)$ is not monitored due to the fact no shortest path between any two vertices of $(S_2\cup \{x\})$ contain edge $e_0$.

        \vspace{0.2cm}

        $\rightarrow$ If $x=(1,\,1)$, then the edge $e_1$ between $(2,2)$ and $(3,2)$ is not monitored. Let $z_1$ $z_2 \in (S_2\cup \{x\})$ monitors edge $e_1$. In this case, $z_1$ or $z_2$ is $(1,1)$ or $(2,2)$. If $z_1=(1,1)$, then the shortest path between $z_1$ and $z_2$ in $G'$ is $P$, where $P=(1,1)  (1,2)(2,2) (3,2) \ldots z_2$. We can find another shortest path, $P'$,  between $z_1$ and $z_2$ in $G'$ from $P$, which does not go via edge $e_1$. Here, $P'=(1,1)(2,1) (3,1)(3,2)\ldots z_2$. Therefore, edge $e_1$ is not monitored by $z_1$ and $z_2$.
        
        Similarly, if $z_1=(2,2)$, then $P=(2,2)(3,2) (3,3) \ldots z_2$ is the shortest path between $z_1$ and $z_2$ in $G'$. We can find another shortest path, $P'$,  between $z_1$ and $z_2$ in $G'$ from $P$, which does not go via edge $e_1$. Here, $P'=(2,2) (2,3) (3,3)\ldots z_2$. Therefore, edge $e_1$ is not monitored by $z_1$ and $z_2$.

\vspace{0.2cm}
        $\rightarrow$ If $x=(2, \,1)$, then the edge $e_2$ between $(2,\,2)$ and $(3,\,2)$ is not monitored. Let $z_1$ $z_2 \in (S_2\cup \{x\})$ monitors edge $e_2$. In this case, $z_1$ or $z_2$ is $(2,1)$ or $(2,2)$. If $z_1=(2,1)$, then the shortest path between $z_1$ and $z_2$ in $G'$ is $P$, where $P=(2,1) (3,1)(3,2) (2,2) (2,3) \ldots  z_2$. We can find another shortest path, $P'$,  between $z_1$ and $z_2$ in $G'$ from $P$, which does not go via edge $e_1$. Here, $P'=(2,1)(1,1) (1,2) (2,2) (2,3) \ldots z_2$. Therefore, edge $e_2$ is not monitored by $z_1$ and $z_2$.
        
        Similarly, if $z_1=(2,2)$, then $P=(2,2)(3,2) (3,3)\ldots z_3$ is the shortest path between $z_1$ and $z_2$ in $G'$. We can find another shortest path, $P'$,  between $z_1$ and $z_2$ in $G'$ from $P$, which does not go via edge $e_1$. Here, $P'=(2,2),(2,3) (3,3) \ldots z_2$. Therefore, edge $e_2$ is not monitored by $z_1$ and $z_2$.

\vspace{0.2cm}
        $\rightarrow$ If $x=(1,\,3)$, then the edge $e_3$ between $(2,\,2)$ and $(1,2)$ is not monitored. Let $z_1$ $z_2 \in (S_2\cup \{x\})$ monitors edge $e_3$. In this case, $z_1$ or $z_2$ is $(3,1)$ or $(2,2)$.
        If $z_1=(3,1)$, then the shortest path between $z_1$ and $z_2$ in $G'$ is $P$, where $P=(3,1) (3,2) (2,2) (1,2)\ldots z_2$. We can find another shortest path, $P'$, between $z_1$ and $z_2$ in $G'$ from $P$, which does not go via edge $e_3$. Here, $P'=(3,1) (2,1) (1,1) (1,2)\ldots z_2$. Therefore, edge $e_3$ is not monitored by $z_1$ and $z_2$.
        
        Similarly, if $z_1=(2,2)$, then $P=(2,2) (1,2)(1,3)\ldots z_2$ is the shortest path between $z_1$ and $z_2$ in $G'$. We can find another shortest path, $P'$,  between $z_1$ and $z_2$ in $G'$ from $P$, which does not go via edge $e_1$. Here, $P'=(2,2)(2,3) (1,3) \ldots z_2$. Therefore, edge $e_3$ is not monitored by $z_1$ and $z_2$.
        
\vspace{0.2cm}
        $\rightarrow$ If $x=(1,2)$ or $(3,2)$, then the edge $e_4$ between $(1,\,1)$ and $(2,\,1)$ is not monitored due to the fact no shortest path between any two vertices of $(S_2\cup \{x\})$ contain edge $e_4$.

         Therefore, we need at least two vertices in set $S_2$ to form a $MEG(G')$, and $S_1$ is the smallest $MEG(G')$. We know $|S_1|=k-2$. Hence, $meg(G')=k-c$.

        \vspace{0.2cm}
        \item \textit{Sub-Case (2)} ${|N(u)\cap C|=1:}$ In this case, $G$ is $m \times n$, where $m,\,n\geq 3$. Without loss of generality, let $u=(2, \,1)$ and $v=(2,2)$. Using Claim \ref{boundary_vertices}, we can show $S-\{(1,1),\,(2, 1),\,(3,1),\,(1,2)\}$ is part of every $MEG(G')$. Similarly, one can show that $v$ is in every $MEG(G')$. We show that $S_1=(S\cup \{v\})-\{(2,1),\,(1,2)\}$ is the smallest $MEG(G')$. Set $S_1$ is forming the $MEG(G')$ due to the following fact. Let us define the shortest path $P_1$, $P_2$, $P_3$, $P_4^l$ and $P_4^{l_1}$ such that $P_1=(1,1)(2,1)(2, 2)$, $P_2=(3,1)(3,2)(2,2)$, $P_3$ is the shortest path between $(2,2)$ and $(2,n)$, $P_4^l$ is the shortest path between $(l,1)$ and $(l, n)$ where $1\leq l\leq m,\, l\neq 2$, and $P_5^{l_1}$ is the shortest path between $(1,l_1)$ and $(m, l_1)$, where $1\leq l_1 \leq n$. It is important to note that these shortest paths are unique in $G'$, and the end vertices of each path are in $S_1$. Each edge $e_1\in G'$ lies on at least one of these unique shortest paths. Therefore, set $S_1$ is a $MEG(G')$. Now, we show that $S_1$ is the smallest $MEG(G')$. The idea for this case is similar to the last case. 

\vspace{0.2cm}
        Note that set $S_2=(S\cup \{v\})-\{(1,1),\,(2, 1),\,(3,1),\,(1,2)\}$ is not a $MEG(G')$ because edge $e$ between $(1,\,1)$ and $(1,\,2)$ is not monitored by set $S_2$. Therefore, we need at least one vertex in set $S_2$. If $(S_2\cup \{ x\})$ forms an $MEG(G')$, then $x \in V(G')-S_2$ (all vertices from $\mathcal{I}$ with $\{(1,1),\,(2, 1),\,(3,1),\,(1,2)\}$).

   \vspace{0.2cm}     
        $\rightarrow$ If $x \in \mathcal{I}$, then the edge $e_0$ between $(1,\,1)$ and $(2,\,1)$ is not monitored due to the fact no shortest path between any two vertices of $(S_2\cup \{x\})$ contain edge $e_0$.
       
\vspace{0.2cm}
        $\rightarrow$ If $x=(1,\,1)$ or $(2,1)$, then the edge $e_1$ between $(3,1)$ and $(3,2)$ is not monitored. Let $x=(1,1)$. Let $z_1$ $z_2 \in (S_2\cup \{x\})$ monitors edge $e_1$. In this case, $z_1$ or $z_2$ is $(1,1)$ or $(i,1)$, $i\geq 4$. If $z_1=(1,1)$, then the shortest path between $z_1$ and $z_2$ in $G'$ is $P$, where $P=(1,1) (2,1), (3,1), (3,2) (3,\,3) \,or\, (4,2) \ldots z_2$. We can find another shortest path, $P'$,  between $z_1$ and $z_2$ in $G'$ from $P$, which does not go via edge $e_1$. Here, $P'=(1,1) (1,2), (2,\,2), (3,2)(3,\,3) \,or\, (4,\,2)  \ldots z_2$. Therefore, edge $e_1$ is not monitored by $z_1$ and $z_2$.
        
        Similarly, if $z_1=(i,1), i\geq 4$, then $P=(i,1)(i-1,1) \ldots(3,1)(3,2)(3,3)\, or \,(4,2)\ldots z_2$ is the shortest path between $z_1$ and $z_2$ in $G'$. We can find another shortest path, $P'$,  between $z_1$ and $z_2$ in $G'$ from $P$, which does not go via edge $e_1$. Here, $P'=(i,1)(i,2) (i-1,2)(3,2)(3,3)\, or \,(4,2)\ldots z_2$. Therefore, edge $e_1$ is not monitored by $z_1$ and $z_2$. A similar argument can be given when $x=(2,1)$. In this case, $z_1$ or $z_2$ is $(2,1)$ or $(i,1)$, $i\geq 4$.

\vspace{0.2cm}
        $\rightarrow$ If $x=(1,2)$, then edge $e_2$ between $(1,1)$ and $(2,1)$ is not monitored. Let $z_1$ $z_2 \in (S_2\cup \{x\})$ monitors edge $e_2$. In this case, the shortest path between $z_1$ and $z_2$ in $G'$ is $P$, where $P=z_1  \ldots (3,1)(2,1)(1,1)(1,2)\ldots z_2$. We can find another shortest path, $P'$, between $z_1$ and $z_2$ in $G'$ from $P$, which does not go via edge $e_2$. Here, $P'=z_1  \ldots (3,1)(3,2)(2,2)(1,2)\ldots z_2$. Therefore, edge $e_2$ is not monitored by $z_1$ and $z_2$.

\vspace{0.2cm}
        $\rightarrow$ If $x=(1,2)$ or $(3,1)$, then edge $e_3$ between $(2,1)$ and $(3,1)$ is not monitored. Let $z_1$ $z_2 \in (S_2\cup \{x\})$ monitors edge $e_2$. Let $x=(1,2)$. In this case, the shortest path between $z_1$ and $z_2$ in $G'$ is $P$, where $P=z_1  \ldots (3,1)(2,1)(1,1)(1,2)\ldots z_2$. We can find another shortest path, $P'$, between $z_1$ and $z_2$ in $G'$ from $P$, which does not go via edge $e_2$. Here, $P'=z_1  \ldots (3,1)(3,2)(2,2)(1,2)\ldots z_2$. Therefore, edge $e_2$ is not monitored by $z_1$ and $z_2$. A similar argument can be given for $x=(3,1)$.

        Therefore, we need at least two vertices in set $S_2$ to form a $MEG(G')$. Therefore, $S_1$ is a smallest $MEG(G')$. We know $|S_1|=k-1$. Therefore, $meg(G')=k-c$.
        
\vspace{0.2cm}
        
       \item \textit{Sub-Case (3)} ${|N(u)\cap C|=0:}$ Without loss of generality, let $u=(i,\,1)$ and $v=(i,2)$, where $2<i<m-1$. Using Claim \ref{boundary_vertices}, we can show that $S-\{(i-1,1), \, (i,1)\, (i+1,1)\}$ is in every $MEG(G')$. Similarly, we can show $v$ is in every $MEG(G')$. It is important to $S_2=(S\cup \{v\})-\{(i-1,1), \, (i,1)\, (i+1,1)\}$ is not a $MEG(G')$ because edge $e_1$ between $(i-1,1)$ and $(i-1,2)$ is not monitored by set $S_2$. Let $P$ be a shortest between $z_1$ and $z_2$ of $S_2$, where $P=z_1 \ldots (i-2,1) (i-1,1)(i-1,2) (i, 2) \ldots z_2$ (or $z_1 \ldots (i, 1)  (i-1,1) (i-1,2)   (i-2,2) \ldots z_2$). In this case, there is another path $P'$ between $z_1$ and $z_2$, which does not go via edge $e_1$. Here, $P'=z_1 \ldots (i-2,1)  (i-2,2) (i-1,2)  (i,2) \ldots z_2$ (or $z_1  \ldots (i, 1) (i,2) (i-1,2)  (i-2,2)\ldots z_2$). Therefore, $S_2$ is not a $MEG(G')$. Therefore, the smallest $MEG(G')$ needs at least one more vertex in $S_2$. We show that $S'=(S\cup \{v\})-\{(i-1,1),\,(i+1,1)\}$ is $MEG(G')$. Let us define the shortest path $P_1$, $P_2$,$P_3$, $P_4^l$, $P_5^{l_1}$ such that $P_1$ is a shortest path between $(i,1)$ and $(i-1,n)$, $P_2$ is the shortest path between $(i,1)$ and $(i+1,n)$, $P_3$ is the shortest path between $(i,2)$ and $(i,n)$, $P_4^l$ is shortest path between $(1,l)$ and $(m,l)$, $1\leq l \leq n$, and $P_5^{l_1}$ is shortest path between $(l_1, 1)$ and $(l_1, n)$, $1\leq l_1\leq m, \, l_1\notin\{i-1,i,i+1\}$. It is important to note that these shortest paths are unique in $G'$, and the end vertices of each path are in $S'$. Each edge $e_1\in G'$ lies on at least one of these unique shortest paths. Therefore, $S'$ is the smallest $MEG(G')$. We know $|S'|=k-1$. Therefore, $meg(G')=k-c-1$.
    \end{itemize}
This completes the proof.
\end{proof}

\begin{lemma}
    Let $G$ be a rectangular $m \times n$ grid graph with $meg(G) = k$, $m, n\geq 3$. Let $G'$ be a partial grid obtained by removing an edge $uv$ from $G$. 
  \begin{equation*}
meg(G') =
\begin{cases} 
k, & \text{if $d_1=d_2=3$,}\\
 k-1, & \text{if $d_1 = 2$, $d_2 = 3$,}\\
   \end{cases}
\end{equation*}
\end{lemma}
\begin{proof}
We divide this proof into two cases as follows.

\noindent \textit{Case (1)} ${d_1=d_2=3}$: Without loss of generality, let $u=(1,j)$ and $v=(1,j+1)$.
Using the argument mentioned in Claim \ref{boundary_vertices}, we can show $S$ contains in every $MEG(G')$. In $G'$, there is a unique shortest path between $(l,1)$ and $(l,n)$ for every $2\leq l \leq m$, there is a unique shortest path between $(1,l_1)$ and $(m,l_1)$ for every $1\leq l_1 \leq m$, there is a unique shortest path between $(1,1)$ and $(1,j)$, and there is a unique shortest path between $(1,j+1)$ and $(1,n)$. If we take the union of all edges of such unique shortest paths, then it is $G'$. Hence, set $S$ is $MEG(G')$, and $S$ is in every $MEG(G')$. Therefore, $S$ is the smallest $MEG(G')$, and $meg(G')=k$.

\noindent \textit{Case (2)} ${d_1=2, \,d_2=3}:$ Without loss of generality, let $u=(1,1)$ and $v=(2,1)$. Note that $m\geq 3$. It is also important to note that $(1,2)$ is not part of the smallest $MEG(G')$ due to $(1,2)$ being a cut vertex. Using a similar argument mentioned in Claim \ref{boundary_vertices}, we can show that $S-\{(1,1),\,(1,2)\}$ is part of every $MEG(G')$, and $(1,1)$ is part of every $MEG(G')$ due to being a pendant vertex. Therefore, $S_1=S-\{(1,2)\}$ is in every $MEG(G')$. Let us define the shortest path $P_1$, $P_2$, $P_3^l$ and $P_4^{l_1}$ such that $P_1$ is shortest path between $(1,1)$ and $ (m,2)$, $P_2$ is shortest path between $(2,1)$ and $(m,1)$, $P_3^l$ is shortest path between $(l,1)$ and $(l, n)$, where $1\leq l\leq m$, $P_4^{l_1}$ is shortest path between $(1,l_1)$ to $(m,l_1)$, where $3 \leq j \leq n$.  It is important to note that these shortest paths are unique in $G'$, and the end vertices of each path are in $S_1$. Each edge $e_1\in G'$ lies on at least one of these unique shortest paths. Therefore, set $S_1$ is a $MEG(G')$, and $S_1$ is the smallest $MEG(G')$, and $meg(G')=k-1$.
This completes the proof.  
\end{proof}

\begin{theorem}
Let $G$ be a rectangular grid of size $m\times n$ with $meg(G)=k$ where $m, n \geq 2$. Let $G'$ be a partial grid obtained by removing an edge $uv$ from $G$. Then, for any $u,v\in V(G)$ such that 
\begin{equation*}
meg(G') =
\begin{cases} 
    k+|S_1|-2, & \text{if $S_1 \neq \emptyset$}\\
    k+2, & \text{if $S_1 =\emptyset$}\\
    k-c, & \text{if $d_1 = 2$, $d_2 = 2$ (or if $d_3=3$ $d_2=4$ and $c\neq 0$}),\\
    k-1, & \text{if $d_3=3,\, d_2=4$ and $c=0$ (or if $d_1 = 2$, $d_2 = 3$, $m,n\geq 3$)}\\
    k, & \text{if $d_1=d_2=3$, $m,n\geq 3$,}\\
\end{cases}
\end{equation*}
where $c$ is either $|N(u)\cap C|$ or $|N(v)\cap C|$ and the set $S_1$ is defined as follows. Let $S_e=\{xy\in E(G): x \in N(u), \, y\in N(v)\; with\; deg(x)=deg(y)=3\}$. If $e=xy\in S_e$, then $x, \,y \in S_1$.
\end{theorem}

\section{Results on General Graphs}\label{general}

Let $G$ be a simple connected graph, and $G'$ be a graph after removing some edge(s) from $G$. In this section, we study how the structural properties of $G$ are affected when specific types of edges are removed, focusing on the changes in the metric $meg(G')$. We will analyze the impact of removing a pendant edge, a cut edge, an edge incident with a cut vertex, and an edge incident with a simplicial vertex. This analysis will highlight how $meg(G')$ varies with different edge removals, shedding light on the graph's structural response to these modifications.

\begin{lemma}\label{lm:upper_pendant}
After removing some (or all) of the pendant edges from a graph \( G \), the resultant graph \( G' \) satisfies:
$$ meg(G') \leq meg(G). $$
\end{lemma}

\begin{proof}
    If $G'$ contains no edge to monitor, then this holds. Let $G'$ contain at least one edge, and $S_G$ be the minimum $MEG(G)$. Let $u$ be a pendant vertex in $G$ and $uv$ be the pendant edge in $G$. To prove our result, it is sufficient to construct one $MEG(G')$ (say $S$), where $G'= G\backslash \{uv\}$ such that $|S| \leq |S_G|$. If $|S| \leq |S_G|$ holds, then $meg(G') \leq meg(G)$ as $|S_{G'}| \leq |S|$, where $S_{G'}$ is the minimum $MEG(G')$. We claim that $S=(S_G \backslash \{u\})\cup \{v\}$ is a $MEG(G')$. Let $e=xy$ be an edge in $G'$. It is important to note that edge $e$ belongs to the edge set of $G$. Let $w$, $w' \in S_G$ be two vertices in $G$ such that $w$ and $w'$ monitor the edge $e$. 

\noindent \textit{Case (1)}  ${ u \notin \{w,\,w'\}}$: In this case, vertices $w$ and $w'$ are in set $S$, and no shortest path $P$ in $G$ between $w$ and $w'$ goes via edge $uv$ due to the fact that if the shortest path $P'$ between $w$ and $w'$ contains edge $uv$, then we can get a smaller path then $P'$. Therefore, every shortest path $P$ between $w$ and $w'$ in $G$ remains the shortest path in $G'$, and vertices of set $S$ monitor $e$.

\noindent \textit{Case (2)} ${ u \in \{w,\,w'\}:}$ Without loss of generality, let $w'=u$. In this case, the edge $e$ is monitored by vertex $u$ and $w$. Since $w$ and $u$ are monitoring edge $e$, therefore it is on all shortest paths between $w$ and $u$. And any shortest path between $w$ and $u$ also goes through vertex $v$ due to $u$ being a pendant vertex of $v$. If there exists a shortest path between $w$ and $v$ in $G'$, which does not contain edge $e$, we can get the shortest path between $w$ and $u$ in $G$, which does not contain edge $e$. Therefore, $e$ is on all shortest paths between $w$ and $v$ in $G'$. Since $w,\,v\in S$, therefore, $e$ is monitored by vertices of set $S$.

    Due to Case $1$, $2$, we can say that $S$ monitors the edges of $G'$ and $|S|=|S_G|$. Since $S_{G'}$ is the minimum $MEG(G')$, $|S_{G'}|\leq |S|$. This completes the proof.   
\end{proof}

\begin{lemma}\label{lm:lower_pendant}
    Let $G$ be a graph with $k$ many pendant edges. If $l$ ($l \leq k$) many pendant edges are removed from $G$, then the following inequality is true.
    $$  meg(G)-l \leq meg(G')$$
\end{lemma}
\begin{proof}
    We will prove it by contradiction. Without loss of generality, let $meg(G')=meg(G)-l-1$ after removing $l$ many pendant edges. Let \( S_1 = \{u_1, u_2, \ldots, u_l\} \) be the pendant vertices, and the edges incident to \( u_i \) for each \( i \in [1, l] \) are removed from \( G \). Let $S_2=\{v_1, \,v_2,\, \ldots , v_{meg(G)-l-1}\}$ 
     be a minimum $MEG(G')$. Let $S$ be $S_1 \cup S_2$. The cardinality of set $S$ is $meg(G)-1$. If we are able to show that $S$ is a $MEG(G)$, then it gives a contradiction. This is because the minimum size of $MEG(G)$ is $meg(G)$. We will show that $S$ is a $MEG(G)$. Let $e=uv$ be an edge of $G$. We have the following two cases.

\noindent \textit{Case (1) ${e\in E(G')}$:} Since $S_2$ is a $MEG(G')$, therefore there exist $v_i$, $v_j\in S_2$ such that the edge $e$ is monitored by $v_i$ and $v_j$ in $G'$. In this case, the edge $e$ is also monitored by $v_i$ and $v_j$ in $G$ due to the following fact. Since \( v_i \) and \( v_j \) are both in \( V(G') \) and the shortest path between them does not pass through any pendant edge incident to \( u_i \), any shortest path \( P \) between \( v_i \) and \( v_j \) in \( G' \) also remains the shortest path in \( G \). The edge $e$ is therefore monitored by the set $S$ in this case.

\noindent \textit{Case (2) ${e\in E(G)-E(G'):}$} In this scenario, edge $e$ is a pendant edge that is removed. In other words, either $u$ or $v$ is $u_i$. Without loss of generality, let $e=uu_i$. Since $v_1 \in S$, edge $e$ lies on all shortest path between $u_i$ and $v_1$. Therefore, edge $e$ is monitored by set $S$. 
   
    Therefore, $S$ is a MEG set of $G$. This completes our proof.
\end{proof}

\begin{theorem}\label{thm:pendant}
    Let $G$ be a graph with $k$ many pendant edges. If $l$ ($l \leq k$) many pendant edges are removed from $G$, then the following inequality is true.
    $$  meg(G)-l \leq meg(G') \leq meg(G).$$
    Moreover, both the lower and the upper bounds are tight.
\end{theorem}
\begin{proof}
    The inequality is true due to Lemma \ref{lm:upper_pendant}, \ref{lm:lower_pendant}. The tightness is as follows.

\medskip

\noindent \textit{Tightness of the upper bound:} Let $G_1$ be a path of length $n\geq 4$. Due to Corollary \ref{cor:existing-path}, $meg(G_1)=2$. After removal of all pendant edges, $meg(G_1')=2$ using Corollary \ref{cor:path}. Using inequality, $0\leq meg(G_1')\leq meg(G_1)$. In $G_1$, $0<meg(G_1')=meg(G_1)$. 

\medskip

 \noindent  \textit{Tightness of the lower bound:} Let $G_2$ be star graph of $n$ vertices. Due to Theorem \ref{thm:existing-Tree}, $meg(G_2)=n-1$. After the removal of all pendant edges ($n-1$ pendant edges), there is no edge to monitor. Therefore, $meg(G_2')=0$. Using inequality, $meg(G_2)-n+1\leq meg(G_2')\leq meg(G_2)$. In graph class $G_2$, $meg(G_2')=meg(G_2)-n+1=0$. 
   
    This completes the proof of the theorem. 
\end{proof}

\begin{figure}
\begin{minipage}{0.45\textwidth}
\centering
\begin{tikzpicture}[scale=0.7,
    level 1/.style={sibling distance=6cm},
    level 2/.style={sibling distance=3cm},
    level 3/.style={sibling distance=1.5cm},
    level distance=1.2cm,
    every node/.style={circle, draw, minimum size=0.6cm, align=center}
]

\node {1}
    child { node {2}
        child { node {4}
            child { node {8} }
            child { node {9} }
        }
        child { node {5}
            child { node {10} }
            child { node {11} }
        }
    }
    child { node {3}
        child { node {6}
            child { node {12} }
            child { node {13} }
        }
        child { node {7}
            child { node {14} }
            child { node {15} }
        }
    };

\end{tikzpicture}
\textbf{(a)}
\end{minipage}
\hfill
\begin{minipage}{0.45\textwidth}
\centering
\begin{tikzpicture}[scale=0.7,
    level 1/.style={sibling distance=6cm},
    level 2/.style={sibling distance=3cm},
    level 3/.style={sibling distance=1.5cm},
    level distance=1.2cm,
    every node/.style={circle, draw, minimum size=0.6cm, align=center}
]

\node {1}
    child { node {2}
        child { node {4}
            child { node {8} edge from parent[draw=none] }
            child { node {9} edge from parent[draw=none] }
        }
        child { node {5}
            child { node {10} edge from parent[draw=none] }
            child { node {11} edge from parent[draw=none] }
        }
    }
    child { node {3}
        child { node {6}
            child { node {12} edge from parent[draw=none] }
            child { node {13} edge from parent[draw=none] }
        }
        child { node {7}
            child { node {14} edge from parent[draw=none] }
            child { node {15} edge from parent[draw=none] }
        }
    };

\end{tikzpicture}
\textbf{(b)}
\end{minipage}

\caption{(a) It is a perfect binary tree for $15$ vertices, (b) It is a perfect binary tree after removal of all pendant edges from the perfect binary of $15$ vertices. }
\label{fig:perfect-Tree}
\end{figure}
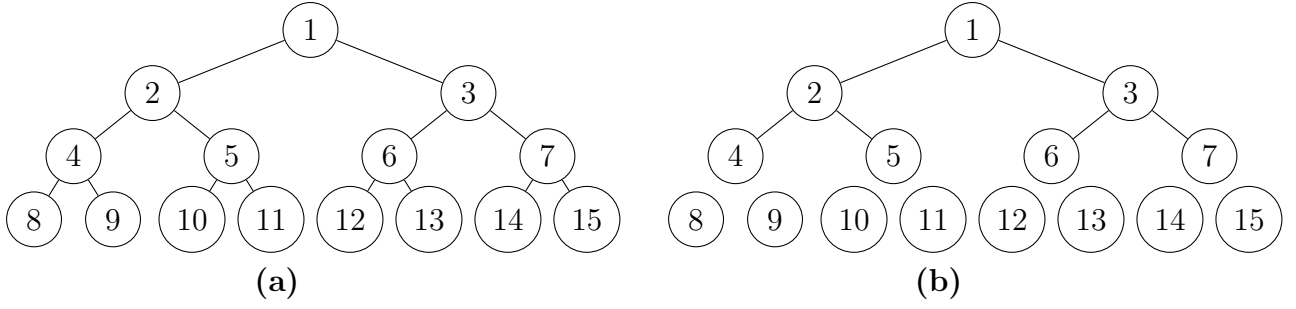

\begin{remark}
    Let $G_3$ be a \textit{perfect binary tree} i.e., a binary tree in which all internal vertices have two children and all leaves are at the same level). Let the size of $4n-1$. We know that if the height of the perfect binary tree is $h$, then the number of vertices is $2^{h+1}-1$, and the number of leaves is $2^h$. If the number of vertices is $4n-1$, then $h=\log(2n)$. Therefore, the number of leaves is $2n$. Using Theorem \ref{thm:existing-Tree}, $meg(G_3)=2n$. The remaining graph is a forest if all the pendant vertices are removed from $G$. Let $G_3'$ be a tree except for all the isolated vertices. It is easy to observe that $G_3'$ is a perfect graph of size $2n-1$. The number of leaves in $G_3'$ is $n$. Due to Theorem \ref{thm:existing-Tree}, $meg(G_3')=n$. In graph class $G_3$, $0<meg(G_3')<2n-1$. One such example for $n=4$ can be seen in Figure \ref{fig:perfect-Tree}.
\end{remark}

\begin{lemma}\label{lm:lower_bridge}
Let $G'=G\backslash \{e\}$, where $e=uv$ is a cut edge in graph $G$ with $u,v\in V(G)$ and $deg(u), deg(v)\geq 2$. The following inequality is true.

$$meg(G) \leq meg(G')$$
\end{lemma}

\begin{proof}
    Let $G_1$ and $G_2$ be two component of $G\backslash \{e\}$. Let $S_1=\{v_1\, \,v_2,\,\ldots,\, v_{k_1}\}$ be a minimum $MEG(G_1)$, and $S_2=\{u_1\, \,u_2,\,\ldots,\, u_{k_2}\}$ be a minimum $MEG(G_2)$. Note that $k_1, \, k_2 \geq 1$ as $deg(u)$, and $deg(v)\geq 2$. As per the definition of minimum monitoring edge-geodetic set of a graph, $S_1 \cup S_2$ is $MEG(G')$. We show that $S_1 \cup S_2$ is a monitoring edge geodetic set of $G$. Let $e_1 \in G$. The following three cases are possible. 

    \noindent \textit{Case (1)} ${e_1 \in G_1:}$ There exists $v_i$ and $v_j$ from set $S_1$ such that $v_i$ and $v_j$ monitor edge $e_1$ in $G_1$ as $S_1$ is a mentoring edge-geodetic set of $G_1$. Since $v_i$, $v_j\in V(G_1)$, the shortest path distance between $v_i$ and $v_j$ in $G_1$ is same as the shortest path distance between $v_i$ and $v_j$ in $G$. Therefore, in this case, any shortest $P$ between $v_i$ and $v_j$ in $G_1$ remains the shortest path in $G$. Also, no shortest path $P_1$ between $v_i$ and $v_j$ exists in $G$, which does not contain edge $e_1$. If such a path $P_1$ exits, then at least one edge of $P_1$ should be edge $uv$; otherwise, it is a path in $G_1$. In this case, $P_1$ can not be the shortest path between $v_i$ and $v_j$ as we find a smaller path between $v_i$ and $v_j$ in $G_1$. Therefore, edge $e_1$ is also monitored by $v_i$ and $v_j$ in $G$.

\noindent \textit{Case (2)} ${e_1 \in G_2:}$ Analogous to Case 1.

\noindent \textit{Case (3)} ${e_1=e:}$ Edge $e$ is monitored by any $v' \in S_1$ and $u' \in S_2$. It is because edge $e$ lies on all shortest paths between $v'$ and $u'$ in $G$.

    \medskip

It implies that $S_1 \cup S_2$ is a $MEG(G)$. It implies that $meg(G) \leq |S_1|+|S_2|=meg(G')$. Therefore, $meg(G) \leq meg(G')$.  This completes the proof.
\end{proof}

\begin{lemma}\label{lm:upper_bridge}
Let $G'=G\backslash \{e\}$, where $e=uv$ is a cut edge in graph $G$ with $u,v\in V(G)$ and $deg(u), deg(v)\geq 2$. The following inequality is true.
$$ meg(G')\leq meg(G)+2$$
\end{lemma}
\begin{proof}
        Let $G_1$ and $G_2$ be two component of $G\backslash \{e\}$, and $S=\{v_1, \, v_2, \, \ldots,\,v_{meg(G)}\}$ be the minimum $MEG(G)$. Let $T_1=(S\cup \{u\})\cap G_1$, and $T_2=(S\cup \{v\})\cap G_2$. We will show that $T_1$ is a $MEG(G_1)$, and $T_2$ is a $MEG(G_2)$. Let $e_1=x_1y_1 \in E(G_1)$, $w, w' \in S$ which are monitoring the edge $e_1$ in $G$. Due to Lemma \ref{lm:existing-cut_vertex}, $u$ and $v$ are not part of set $S$. Therefore, $w$ and $w'$ can't be $u$ or $v$.

        \noindent \textit{Case (1)} ${w,\,w'\in V(G_1):}$ Any shortest path between $w$ and $w'$ in $G$ is also the shortest path in $G_1$, because if any path $P$ between $w$ and $w'$ is not in $G_1$, then it contains the edge $uv$, and we can get a smaller path between $w$ and $w'$. Therefore, the edge $e_1$ is monitored by the set $T_1$.

         \noindent \textit{Case (2)} ${w \in V(G_1)}$, ${w' \in V(G_2):}$ All shortest paths between $w$ and $w'$ in $G$ contain edge $uv$. Therefore, edge $e$ also lies on all shortest paths between $w$ and $u$ in $G_1$. Since $w, \, u \in T_1$, therefore $e_1$ is monitored by set $T_1$.

         \noindent  \textit{Case (3)} ${w' \in V(G_1)}$, ${w \in V(G_2):}$ Analogous to Case 2.

     Similarly, we can show that $T_2$ is $MEG(G_2)$. Since, $meg(G')\leq |T_1|+|T_2|=meg(G)+2$, therefore, $meg(G')\leq meg(G)+2$. This completes our proof.
\end{proof}

Using Lemma \ref{lm:lower_bridge} and Lemma \ref{lm:upper_bridge}, we have the following theorem.

\begin{theorem}\label{thm:bridge}
    Let $G'=G\backslash \{e\}$, where $e=uv$ is a cut edge in graph $G$ with $u,v\in V(G)$ and $deg(u), deg(v)\geq 2$. The following inequality holds.
$$ meg(G)\leq meg(G')\leq meg(G)+2.$$ 
Moreover, both the lower and the upper bounds are tight.
\end{theorem}
\begin{proof}
    The inequality is true due to Lemma \ref{lm:upper_bridge}, \ref{lm:lower_bridge}. The tightness is as follows. 

   \medskip
   
    \noindent \textit{Tightness of the upper bound:} Let $H_1$ be a tree of size $n>6$ which has at least one edge $e=uv$ such that $deg(u)$ and $deg(v)$ is 2. Using Theorem~\ref{lm:tree-1-dyn}, if the edge $e$ is removed, then $meg(H_1')=meg(H_1)+2$. Therefore, graph class $H_1$ shows the tightness of the upper bound of the mentioned inequality.
       
\medskip
        
    \noindent \textit{Tightness of the lower bound:} Let $H_2$ be a tree of size $n>6$ which has at least one edge $e=uv$ such that $deg(u)$ and $deg(v)$ is at least 3. Using Theorem \ref{lm:tree-1-dyn}, if the edge $e$ is removed, then $meg(H_2)=meg(H_2')$. Therefore, graph class $H_2$ shows the tightness of the lower bound of the mentioned inequality.

    This completes the proof of the theorem.
    \end{proof}

\begin{remark}
    Let $H_3$ be a tree of size $n>6$ which has at least one edge $e=uv$ such that $deg(u)$ is at least 3, and $deg(v)$ is 2. Using Theorem \ref{lm:tree-1-dyn}, if the adversary removes edge $e$, then $meg(H_3')=meg(H_3)+1$. Therefore, graph class $H_3$ satisfies $meg(H_3)<meg(H_3')<meg(H_3)+2$.
\end{remark}

In \cite{foucaud2023monitoringfull}, the authors discuss the MEG-extremal graphs. A graph $G$ is called MEG-extremal if $meg(G)=n$, where $n$ is the number of vertices in $G$. The following question arises for MEG-extremal graphs.

\vspace{0.2cm}

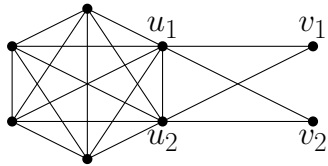
\begin{figure}[h]
    \centering
    \begin{tikzpicture}[scale=0.5, transform shape]
    \node[noeud,fill=black,circle] (2) at (-13, 2){};
     \node[noeud,fill=black,circle] (3) at (-11, 5){};
    \node[noeud,fill=black,circle] (4) at (-11, 1){};
    \node[noeud,fill=black,circle] (5) at (-13, 4) {};
    \node[noeud,fill=black,circle] (6) at (-9, 4) {};
    \node at (-9, 4.5) {\huge $u_1$};
    \node[noeud,fill=black,circle] (7) at (-9, 2) {};
    \node at  (-9, 1.5) {\huge $u_2$};
    \node[noeud,fill=black,circle] (8) at (-5, 4) {};
     \node at  (-5, 4.5) {\huge $v_1$};
    \node[noeud,fill=black,circle] (9) at (-5, 2) {};
     \node at  (-5, 1.5) {\huge $v_2$};
        
        \draw (6) -- (9);
        \draw (7) -- (8);
        \draw (5) -- (3);
        \draw (6) --(8);
        \draw (7)--(9);
        \draw (5) -- (2);
        \draw (2) -- (4);
        \draw (4) -- (3);
        \draw (3) -- (6);
        \draw (6) -- (7);
        \draw (7) -- (4);
        \draw (2) -- (6);
        \draw (5) -- (7);
        \draw (2) -- (3);
        \draw (2) -- (7);
        \draw (3) -- (7);
        \draw (6) -- (4);
        \draw (5) -- (4);
        \draw (5) -- (6);
    \end{tikzpicture}
    \caption{Construction of split graph $G$ for $n=8$.}
    \label{fig:split}
\end{figure}

\begin{question}\label{answer}
If an edge $e$ is removed from the MEG-extremal graph $G$ with order $n$, does it follow that $G - \{ e \}$ remains a MEG-extremal graph?
\end{question}
\begin{answer}
    The answer to this question is negative. We will construct a graph, denoted as $G$, which is the MEG-extremal graph. In this graph $G$, there exists an edge $e$ such that if the edge $e$ is removed from $G$ with order $n$, then $G\backslash \{e\}$ is no longer a MEG-extremal graph.

    Let $K_{n-2}$ be a clique of size $n-2$, and let $u_1$ and $u_2$ be two vertices of $K_{n-2}$. We introduce two new vertices, $v_1$ and $v_2$. In the graph $G$, $v_1$ is adjacent to $u_1$, $u_2$ of $K_{n-2}$, and $v_2$ is also adjacent to $u_1$, $u_2$ of $K_{n-2}$ (refer to Fig. \ref{fig:split} for $n=8$). It is not difficult to observe that $G$ is a split graph. Using Corollary \ref{cor:exiting-split}, $meg(G)=n$. If any edge from set $\{u_1v_1,\, u_1v_2, \, u_2v_1,\, u_2v_2, \}$ is removed, then there will be a vertex with pendant neighbour. In this case, using Corollary \ref{cor:exiting-split}, $meg(G')<n$. Therefore, if an edge $e$ is removed from the MEG-extremal graph, say $G$, then $G\backslash \{e\}$ is not necessarily a MEG-extremal graph.  
\end{answer}

\begin{theorem}
   Let $G$ be a graph and $v$ be a cut vertex in $G$. The following inequality holds after removing an edge $e$ incident to vertex $v$.
   
   $$ meg(G')\leq meg(G)+2$$
   Here, $G'=G\backslash \{e\}$, and moreover, this bound is tight.
\end{theorem}
\begin{proof}
    Let $e=uv$ be an edge in $G$ which is incident to cut vertex $v$. If $e$ is a cut edge, then due to Theorem \ref{thm:bridge}, $meg(G') \leq meg(G)+2$. Therefore, $meg(G') \leq meg(G)+2$. Suppose $e$ is not a cut edge in $G$. Let $S$ be a minimum $MEG(G)$. Let $S'=S\cup \{u, v\}$. If $S'$ is a $MEG(G')$, then $meg(G') \leq |S'|=meg(G)+2$. Therefore, $meg(G') \leq meg(G)+2$. We show that $S'$ is a $MEG(G')$. After removing an edge $e$ from $G$, the shortest path distance between some vertices can be increased. Let $e'$ be an edge in $G'$. Edge $e'$ is in $E(G')$ due to the fact that $G'$ is a sub-graph of $G$. Since $S$ is $MEG(G)$, there exist $x_1,\,x_2\in S$ such that edge $e'$ is monitored by $x_1$ and $x_2$ in $G$.

    If after removing edge $e$ from $G$, the shortest distance between $x_1$, $x_2$ does not increase, then $e'$ is monitored by $x_1$ and $x_2$ in $G'$. Since $x_1,\,x_2\in S'$, $e'$ is monitored by set $S'$. If after removing edge $e$, the shortest path distance between $x_1$ and $x_2$ is increased, then the edge $e=uv$ is in all the shortest path between $x_1$ and $x_2$ in $G$. Therefore, all shortest paths $P$ between $x_1$ and $x_2$ look like $x_1\ldots uv \ldots x_2$. Edge $e'$ either is in $P_1=x_1 \ldots u$ or $P_2=v \ldots x_2$. Note that the shortest path distance between $u$ and $x_1$ in $G'$ is the same as $P_1$, and the shortest path distance between $v$ and $x_2$ in $G'$ is the same as $P_2$. If not, we can find a path in $G$ whose length is smaller than $P$. Let $P_1=x_1 \ldots u$ be the shortest path between $x_1$ and $u$ in $G'$, which does not contain edge $e'$, and $P_2=x_2 \ldots v$ be the shortest path between $x_2$ and $v$ in $G'$ which does not contain edge $e'$. Let $P=x_1 \ldots uv\ldots x_2$ be a path in $G$ which does not contain edge $e'$. The path $P$ is the shortest path in $G$ due to the fact that the shortest path distance between $u$ and $x_1$ is the same as $P_1$, and the shortest path distance between $v$ and $x_2$ is the same as $P_2$. Therefore, the shortest path $P$ between $x_1$ and $x_2$ in $G$ exists, which does not contain edge $e'$. This is a contradiction. Therefore, if $e'$ is in $P_1$, then edge $e'$ is monitored by $x_1$ and $u$ in $G'$; for the other case, $e'$ is monitored by $v$ and $x_2$ in $G'$. Since $u$, $v \in S'$. It implies that $e'$ is monitored by $x_1$ and $u$ (or $x_2$ and $v$). Therefore, $S'$ is a $MEG(G')$. This completes the proof. 

    \medskip
    
    \textit{Tightness of the bound:} Let $G$ be a tree, and there exist an edge $uv$ such that $deg(u)=deg(v)=2$. In this case, due to Theorem \ref{lm:tree-1-dyn}, $meg(G')=meg(G)+2$.
    
\end{proof}

\begin{lemma}\label{lm:simplicial-upper}
    Let $G$ be a graph, and let $v$ be a simplicial vertex in $G$ with $deg(v)\geq 2$. The following inequality is true after removing an edge incident to vertex $v$.
    $$meg(G')\leq meg(G)+1$$
\end{lemma}
\begin{proof}
    Due to Lemma \ref{lm:existing-simplicial}, $v$ is part of every $MEG(G)$. Let $G'=G\backslash \{e_1=uv\}$, where $e_1$ is an incident edge of $v$ in $G$, and $M=\{u_1, u_2\ldots u_{deg(v)}\}$ be a set of neighbours of $v$ in $G$. It is important to note that after removing edge $e$ incident to $v$ from $G$, $v$ remains a simplicial vertex. Therefore, it is part of every $MEG(G')$. Let $S$ be the smallest $MEG(G)$. If we show that $S_1=S\cup \{u\}$ is a $MEG(G')$, then $meg(G')\leq |S_1|=|S\cup \{u\}|\leq meg(G)+1$. We show that $S_1$ is a $MEG(G')$. 
    
    Let $e$ be an edge in $G'$. Since $G'$ be a sub-graph of $G$, $e$ is an edge of $G$. Therefore, there exist $w$ and $w'$ in $S$, such that edge $e$ is monitored by $w$ and $w'$ in $G$, i.e., $e$ lies on all shortest path $P$ between $w$ and $w'$. Let $P=(w=)x_1 x_2 x_3 \ldots x_k(=w')$ be a shortest path between $w$ and $w'$ in $G$.

    \noindent \textit{Case (1)} ${w\neq v}$ \textit{and} ${w'\neq v:}$ In this case, $x_i \neq v$ due to the following fact: If $x_i=v$, then $x_{i-1},\, x_{i+1}\in M$. Since $v$ is a simplicial vertex, there exists an edge between $x_{i-1}$ and $x_{i+1}$. In this case, $P$ can not be the shortest path between $w$ and $w'$ due to the fact that there exists a path $P'=x_1 x_2 \ldots x_{i-1} x_{i+1}\ldots x_k$ between $w$ and $w'$, and its length is $k-1$. Therefore, if $w\neq v$ and $w'\neq v$, then $x_i \neq v$. It means that any shortest path $P$ in $G$ remains the shortest path in $G'$. Therefore, in this case, edge $e$ is monitored in $G'$ as well due to $w$, $w'\in S_1$.

    \noindent \textit{Case (2)} ${w=v}$ \textit{or} ${w'= v:}$ Without loss of generality, let $w'=v$. In this case, $x_k$ is $v$ in any shortest path $P$. After the removal of edge $uv$, the shortest path distance between $w$ and $v$ can be increased. If the shortest path distance between $w$ and $v$ does not increase, then edge $e$ is monitored in $G'$ as well due to $w,\, v \in S_1$. If the shortest path distance between $w$ and $v$ increases, then $x_{k-1}$ is $u$ for every shortest path $P$ between $w$ and $v$. In this case, the shortest path distance between $w$ and $u$ in $G'$ can not be less than $k-1$; otherwise, a path exists between $w$ and $v$ whose length is less than path $P$. Therefore, the shortest path distance between $w$ and $u$ in $G'$ is $k-1$. In this case, a path $P'=(w=)x_1 x_2 x_3 \ldots x_{k-1}(=u)$ be a shortest path between $w$ and $u$ as well, where $P'$ is a sub-path of $P$, and edge $e$ lies on path $P'$. We claim that all such sub-path $P'$ of $P$ are the only shortest path between $w$ and $u$ in $G'$. If there exists another shortest path $P''$ in $G'$ such that edge $e$ does not lie on $P''$, then $w$ and $v$ can not monitor edge $e$ in $G$ due to the fact that adding vertex $v$ at the end of path $P''$ gives the shortest path of length $k$ in $G$ which does not contain edge $e$. Therefore, $P'$ is the shortest path and edge $e$ lies on path $P'$. Since $w, u \in S_1$. Therefore, edge $e$ is monitored by $w$ and $u$. This completes the proof. 
\end{proof}

\begin{lemma}\label{lm:simplicial-lower}
    Let $G$ be a graph and $v$ be a simplicial vertex in $G$ with $deg(v)\geq 2$. The following inequality is true after removing an edge incident to vertex $v$.
    $$meg(G)-deg(v) \leq meg(G')$$
\end{lemma}
\begin{proof}
    Let $G'=G\backslash \{e_1=uv\}$, where $e_1$ is an incident edge of $v$ in $G$, and $S$ be a smallest $MEG(G')$, and $M=\{u_1, u_2\ldots u_{deg(v)}\}$ be a set of neighbours of $v$ in $G$. If we show that $S_1=S\cup M$ be a $MEG(G)$, then $meg(G)\leq |S\cup M|=meg(G)+deg(u)\implies meg(G)-deg(u)\leq meg(G')$. Note that it may be possible that $u_i\in S$ for some $i$. It is also important to note that after removing edge $e$ incident to $v$ from $G$, $v$ remains a simplicial vertex. Therefore, $v \in S$ using Lemma \ref{lm:existing-simplicial}. Let $e$ be an edge of $G$. We divide this proof into two parts.

         \noindent\textit{Case (1)} ${e=uv:}$ Since $u$ and $v\in S_1$, and there is path of length 1 in $G$, therefore, edge $e$ is monitored by set $S_1$ in $G$.

        \noindent \textit{Case (2)} ${e\neq uv:}$ In this case, $e\in G'$. Therefore, $w$ and $w'$ exist in set $S$ such that edge $e$ lies on all shortest path $P$ between $w$ and $w'$ in $G'$. There are two sub-cases.

       \noindent \textit{Sub-Case (1)} ${w\neq v}$ and ${w'\neq v:}$ 
            
            \noindent ${\rightarrow}$ If $w$ and $w' \in M$, then $e=ww'$. In this case, edge $e$ is monitored by set $S_1$. 

            \noindent ${\rightarrow}$ If $w' \in M$ and $P$ is the shortest path in $G'$, then $w, \,w'$ also monitors in $G$. If not, then there exists a path $P'$ between via edge $uv$, which does not contain the edge $e$. In this case, $P'=x_1(=w) x_2 \ldots u v w'$. It may be possible that $w'=u$. In this case, there exists one path $P''=x_1(=w) x_2 \ldots u$ whose length is smaller than $P'$. Therefore, if $w'=u$, then the existence of path $P'$ is not possible. If $w'\neq u$, then there exist one path $P'''=x_1(=w) x_2 \ldots uw'$ whose length is less than length of $P'$ due the fact $u,\, w'\in M$. Therefore, if $w'\neq u$, such a path $P'$ does not exist. Similarly, we can show for $w \in M$ and $w'\notin M$. 

            \noindent ${\rightarrow}$ If $w, \, w' \notin M$ and $P$ is the shortest path in $G'$, then $w, \,w'$ also monitors in $G$. If not, then there exists a path $P_1$ between via edge $uv$, which does not contain the edge. In this case, $P_1=x_1(=w) x_2 \ldots u v u_i \ldots w'$, where $u_i\in M$. Since $u,\, u_i \in M$, therefore there exist one path $P_2=x_1(=w) x_2 \ldots uw'$ whose length is less than length of $P_1$. Therefore, such a path $P_1$ does not exist.

             \noindent \textit{Sub-Case (2)} ${w=v}$ or ${w'= v:}$ Without loss of generality, let $w'=v$. Let $k$ be the length of path $P$. Since, edge $e$ is monitored by $w$ and $v$ in $G'$, edge $e$ lies on all $P$ between $w$ and $v$ in $G'$. In this case, $P=y_1(=w) y_2 \ldots u_i v$, where $u_i \in M$. Note that $u_i$ can not be $u$ as $uv \notin E(G')$. Edge $e$ lies on all shortest paths between $w$ and $u_i$ in $G'$. If not, we can find a path between $w$ and $v$ in $G'$ of length $k$, which does not contain edge $e$, which is impossible. Edge $e$ is also monitored by $w$ and $u_i$ in $G$. If not, there exists a shortest path $P'$ between $w$ and $u_i$ in $G$ such that it does not contain edge $e$. Let $P'=x_1(=w) x_2 \ldots u_i$. Note that $P'$ does not contain edge $uv$. If it does, then $P'=x_1(=w) x_2 \ldots u v u_i$. Since $u_i, \,u\in M$, there exists a path between $w$ and $u_i$ with less than the length of path $P'$. Therefore, $P'$ does not contain edge $uv$. If $P'$ is the shortest path in $G$, which does not contain edge $uv$, then $P'$ is also the shortest path $G'$. In this case, adding a vertex at the end of path $P'$ gives the shortest path between $w$ and $v$ in $G'$, which does not contain edge $e$. This gives a contradiction. Therefore, edge $e$ is monitored by $w$ and $u_i$ in $G$. This completes the proof. 

\end{proof}

We have the following theorem based on Lemma \ref{lm:simplicial-upper}, \ref{lm:simplicial-lower}.
\begin{theorem}
    Let $G$ be a graph, and let $v$ be a simplicial vertex in $G$ with $deg(v)\geq 2$. If an edge $e$ incident to vertex $v$ is removed, then the following inequality is true. $$meg(G)-deg(v)\leq meg(G')\leq meg(G)+1$$ 
\end{theorem}

\section{Conclusion}\label{sec:conclusion}

In this work, we studied the impact of edge removals on the monitoring edge-geodetic number, $meg(G)$, across different graph classes. We establish bounds on $meg(G)$ after the removal of specific types of edge(s) from general graphs, including pendant edges, cut edges, the edges incident to cut vertices, and the edges corresponding to simplicial vertices.  Our results reveal that structural changes resulting from edge deletions can significantly alter the composition and the minimum size of the $MEG$-set.

Moreover, This work provides insights into the behaviour of $meg(G)$ under edge removal and opens up several directions for future research. Future studies could analyze different dynamic graph types: the impacts of vertex removal, or the addition of edges. Additionally, the development of efficient algorithms for monitoring edge-geodetic sets in dynamic settings would greatly improve the resilience and adaptability of network designs.

\bibliography{reference}

\begin{thebibliography}{10}

\bibitem{Aaron_2011}
Aaron Bernstein and Liam Roditty.
\newblock Improved dynamic algorithms for maintaining approximate shortest paths under deletions.
\newblock In {\em Proceedings of the Twenty-Second Annual ACM-SIAM Symposium on Discrete Algorithms}, SODA '11, page 1355–1365, USA, 2011. Society for Industrial and Applied Mathematics.

\bibitem{Sayan_2017}
Sayan Bhattacharya, Deeparnab Chakrabarty, and Monika Henzinger.
\newblock Deterministic fully dynamic approximate vertex cover and fractional matching in o(1) amortized update time.
\newblock In Friedrich Eisenbrand and Jochen Koenemann, editors, {\em Integer Programming and Combinatorial Optimization}, pages 86--98, Cham, 2017. Springer International Publishing.

\bibitem{Davide_2024}
Davide Bilò, Giordano Colli, Luca Forlizzi, and Stefano Leucci.
\newblock On the inapproximability of finding minimum monitoring edge-geodetic sets, 2024.
\newblock URL: \url{https://arxiv.org/abs/2405.13875}, \href {https://arxiv.org/abs/2405.13875} {\path{arXiv:2405.13875}}.

\bibitem{Camil_2004}
Camil Demetrescu and Giuseppe~F. Italiano.
\newblock A new approach to dynamic all pairs shortest paths.
\newblock {\em J. ACM}, 51(6):968–992, November 2004.
\newblock \href {https://doi.org/10.1145/1039488.1039492} {\path{doi:10.1145/1039488.1039492}}.

\bibitem{foucaud2023monitoringfull}
Subhadeep~R. Dev, Sanjana Dey, Florent Foucaud, Narayanan Krishna, and Lekshmi~Ramasubramony Sulochana.
\newblock Monitoring edge-geodetic sets in graphs, 2023.
\newblock arXiv preprint 2210.03774.
\newblock URL: \url{https://arxiv.org/abs/2210.03774}, \href {https://arxiv.org/abs/2210.03774} {\path{arXiv:2210.03774}}.

\bibitem{Foucad_2024}
Florent Foucaud, Pierre-Marie Marcille, Zin~Mar Myint, R.~B. Sandeep, Sagnik Sen, and S.~Taruni.
\newblock Monitoring edge-geodetic sets in graphs: Extremal graphs, bounds, complexity.
\newblock In {\em Algorithms and Discrete Applied Mathematics: 10th International Conference, CALDAM 2024, Bhilai, India, February 15–17, 2024, Proceedings}, page 29–43, Berlin, Heidelberg, 2024. Springer-Verlag.
\newblock \href {https://doi.org/10.1007/978-3-031-52213-0_3} {\path{doi:10.1007/978-3-031-52213-0_3}}.

\bibitem{foucaud2024bounds}
Florent Foucaud, Pierre-Marie Marcille, Zin~Mar Myint, RB~Sandeep, Sagnik Sen, and S~Taruni.
\newblock Bounds and extremal graphs for monitoring edge-geodetic sets in graphs.
\newblock {\em arXiv preprint arXiv:2403.09122}, 2024.

\bibitem{foucaud2023monitoring}
Florent Foucaud, Krishna Narayanan, and Lekshmi Ramasubramony~Sulochana.
\newblock Monitoring edge-geodetic sets in graphs.
\newblock In {\em Algorithms and Discrete Applied Mathematics: 9th International Conference, CALDAM 2023, Gandhinagar, India, February 9--11, 2023, Proceedings}, pages 245--256. Springer, 2023.

\bibitem{harary1997dynamic}
Frank Harary and Gopal Gupta.
\newblock Dynamic graph models.
\newblock {\em Mathematical and Computer Modelling}, 25(7):79--87, 1997.

\bibitem{haslegrave2023monitoring}
John Haslegrave.
\newblock Monitoring edge-geodetic sets: hardness and graph products.
\newblock {\em Discrete Applied Mathematics}, 340:79--84, 2023.

\bibitem{Jacob_1998}
Jacob Holm, Kristian de~Lichtenberg, and Mikkel Thorup.
\newblock Poly-logarithmic deterministic fully-dynamic algorithms for connectivity, minimum spanning tree, 2-edge, and biconnectivity.
\newblock In {\em Proceedings of the Thirtieth Annual ACM Symposium on Theory of Computing}, STOC '98, page 79–89, New York, NY, USA, 1998. Association for Computing Machinery.
\newblock \href {https://doi.org/10.1145/276698.276715} {\path{doi:10.1145/276698.276715}}.

\bibitem{Li_2024}
Xin Li, Wen Li, Ao~Tan, Mengmeng He, and Weizhen Chen.
\newblock Monitoring edge-geodetic numbers of mycielskian graph classes.
\newblock {\em Journal of Interconnection Networks}, 0(0):2450010, 0.
\newblock \href {https://doi.org/10.1142/S0219265924500105} {\path{doi:10.1142/S0219265924500105}}.

\bibitem{Ma_2024}
Yifan~Yao Rongrong~Ma, Zhen~Ji and Yalong Lei.
\newblock Monitoring-edge-geodetic numbers of radix triangular mesh and sierpiński graphs.
\newblock {\em International Journal of Parallel, Emergent and Distributed Systems}, 39(3):353--361, 2024.
\newblock \href {https://doi.org/10.1080/17445760.2023.2294369} {\path{doi:10.1080/17445760.2023.2294369}}.

\bibitem{Sankowski_2004}
Piotr Sankowski.
\newblock Dynamic transitive closure via dynamic matrix inverse: extended abstract.
\newblock {\em 45th Annual IEEE Symposium on Foundations of Computer Science}, pages 509--517, 2004.
\newblock URL: \url{https://api.semanticscholar.org/CorpusID:14447500}.

\bibitem{xu2024monitoring}
Xin Xu, Chenxu Yang, Gemaji Bao, Ayun Zhang, and Xuan Shao.
\newblock Monitoring-edge-geodetic sets in product networks.
\newblock {\em International Journal of Parallel, Emergent and Distributed Systems}, 39(2):264--277, 2024.
\newblock \href {https://doi.org/10.1080/17445760.2024.2301929} {\path{doi:10.1080/17445760.2024.2301929}}.

\end{thebibliography}
\end{document}